\documentclass[review]{elsarticle}

\usepackage{lineno,hyperref,comment,xcolor}\def\d{{\, \rm d}}
\usepackage{amssymb,amsthm,amsmath,bm}\newtheorem{prop}{Proposition}
\usepackage{soul}
\usepackage{epstopdf}
\modulolinenumbers[5]

\journal{Physica D: Nonlinear Phenomena}

\begin{document}

\begin{frontmatter}

\title{A Causality-Based Learning Approach for Discovering the Underlying Dynamics of Complex Systems from Partial Observations with Stochastic Parameterization }

\author{Nan Chen}
\address{Department of Mathematics, University of Wisconsin-Madison, 480 Lincoln Dr., Madison, WI 53706, USA}
\author{Yinling Zhang\corref{mycorrespondingauthor}}\cortext[mycorrespondingauthor]{Corresponding author}\ead{zhang2447@wisc.edu}
\address{Department of Mathematics, University of Wisconsin-Madison, 480 Lincoln Dr., Madison, WI 53706, USA}





\begin{abstract}
Discovering the underlying dynamics of complex systems from data is an important practical topic. Constrained optimization algorithms are widely utilized and lead to many successes. Yet, such purely data-driven methods may bring about incorrect physics in the presence of random noise and cannot easily handle the situation with incomplete data. In this paper, a new iterative learning algorithm for complex turbulent systems with partial observations is developed that alternates between identifying model structures, recovering unobserved variables, and estimating parameters. First, a causality-based learning approach is utilized for the sparse identification of model structures, which takes into account certain physics knowledge that is pre-learned from data. It has unique advantages in coping with indirect coupling between features and is robust to the stochastic noise. A practical algorithm is designed to facilitate the causal inference for high-dimensional systems. Next, a systematic nonlinear stochastic parameterization is built to characterize the time evolution of the unobserved variables.
Closed analytic formula via an efficient nonlinear data assimilation is exploited to sample the trajectories of the unobserved variables, which are then treated as synthetic observations to advance a rapid parameter estimation. Furthermore, the localization of the state variable dependence and the physics constraints are incorporated into the learning procedure, which mitigate the curse of dimensionality and prevent the finite time blow-up issue. Numerical experiments show that the new algorithm succeeds in identifying the model structure and providing suitable stochastic parameterizations for many complex nonlinear systems with chaotic dynamics, spatiotemporal multiscale structures, intermittency, and extreme events.
\end{abstract}


\end{frontmatter}

\section{Introduction}
Complex turbulent dynamical systems appear in many areas, such as geophysics, neural science, engineering, and atmosphere and ocean science \cite{strogatz2018nonlinear, sheard2009principles, ghil2012topics}. These complex systems are characterized by strong nonlinearity, high dimensionality, and multiscale structures. The nonlinear interactions between different scales transfer energy throughout the system, which triggers a large dimension of strong instabilities. As a result, many non-Gaussian features, such as extreme and rare events, intermittency and fat-tailed probability density functions (PDFs), are observed in these systems \cite{sapsis2021statistics, lucarini2016extremes, franzke2015stochastic, wilcox1988multiscale, majda2016introduction, tao2009multiscale}. In addition to improving the description of the phenomena, accurate modeling of these complex systems is the prerequisite of state estimation, uncertainty quantification, and prediction \cite{kalnay2003atmospheric, lahoz2010data, majda2012filtering, evensen2009data, law2015data}.

Due to the incomplete physical understanding of nature and the inadequate model resolution resulting from the limited computing power, model error is often inevitable in inferring these complex systems utilizing the purely knowledge-based modeling approaches that aim at revealing the entire model structure based on explicit physical laws \cite{palmer2001nonlinear, majda2012lessons, orrell2001model, hu2010ensemble, benner2015survey}. To this end, learning the underlying dynamics of these complex systems with the help of data is of practical significance. In fact, exploiting the data in an appropriate fashion may facilitate the discovery of additional physical structures beyond those obtained from the available partial knowledge. It also provides useful information for the development of effective parameterizations that allows to compensate for the inadequate model resolution.
Depending on the amount of available physical understanding of the problem of interest, the data-driven learning methods can be divided into two categories.
On one hand, if there is little prior knowledge of the underlying dynamics, then the model structure and the model parameters have to be learned almost completely from data. In such a situation, the learning algorithm typically starts with a large library of candidate functions, and then a certain sparse identification technique, such as the least absolute shrinkage and selection operator (LASSO) regression \cite{santosa1986linear, tibshirani1996regression}, is incorporated into the optimization procedure for the function selection to obtain a parsimonious model \cite{brunton2016discovering, boninsegna2018sparse, schneider2021learning, cortiella2021sparse, schaeffer2018extracting, goyal2021learning}.
On the other hand, when part of the underlying dynamics is known, then the task becomes more straightforward as the learning algorithm is primarily utilized to discover the residual part. In general, simple closure terms or parameterizations are developed from data to approximate the residual such that the complexity of the discovered model does not increase significantly \cite{ahmed2021closures, chekroun2017data, chekroun2021stochastic, xie2020closure, chen2018conditional, lin2021data, mou2021data, peherstorfer2015dynamic, hijazi2020data, smarra2018data,  majda2012physics, harlim2014ensemble, chen2022shock, mou2022efficient, kondrashov2015data, chattopadhyay2020superparameterization, majda2018model}. Note that if the available partial dynamics are inaccurate, then systematic learning algorithms can be developed to either correct the model error explicitly from data \cite{mojgani2021closed} or introduce additional judicious model error to offset the existing bias \cite{gershgorin2010improving, gottwald2013role, branicki2013non}. It is also worthwhile to mention that sometimes the learning output is represented in the non-parametric forms, such as the neural networks, if the main goal is to forecast the system instead of reaching the explicit physical formulation \cite{wouters2013multi, santos2021reduced, san2018extreme, chattopadhyay2020superparameterization, chen2021bamcafe, pawar2020data, moosavi2015efficient, chen2022physics, chattopadhyay2020deep, chattopadhyay2021towards}.

These data-driven methods have led to many successes in various contexts, especially in building approximate models and forecasting time series. Yet, several challenges still exist in exploiting data-driven approaches to robustly discover the underlying physics of complex turbulent systems. First, satisfying the model parsimony is only a necessary but not sufficient condition for identifying the true dynamics. In fact, the model identification based on the purely data-driven constrained optimization algorithms does not take into account the explicit physics. As a result, in the presence of even small random noise, both the covariate selection accuracy and the fraction of zero entries may decrease significantly  \cite{elinger2020information}, leading to a large bias in discovering the underlying physics. Second, it is often the case that only the observations of a subset of the state variables are available in practice, which is known as the partial observations. In such a situation, the state estimation of the unobserved variables and the parameterization of these processes have to be carried out at the same time as the discovery of the underlying dynamics of the observed variables and the parameter estimation. This leads to a significant increase of the computational cost since the uncertainty quantification of the estimated unobserved states has to be incorporated into the optimization procedure. In addition, as the dimension of the system becomes large, the number of candidate functions in the library often shoots up as well. The curse of dimensionality prevents an efficient selection of the most relevant functions.

In this paper, a causality-based iterative learning algorithm is developed, which aims at overcoming the above difficulties in discovering the dynamics of complex turbulent systems from only partial observations. The algorithm alternates between identifying model structures, recovering unobserved variables, and estimating parameters.
First, the model identification procedure is different from the LASSO regression and many other straightforward constrained optimizations, where data is directly utilized to compute the loss function that involves a regularizer for the model sparsity. In the proposed approach, a causality-based sparse identification of the model structure is adopted, which takes into account certain physics knowledge that is pre-learned from data. Specifically, in light of the observational data, an information measurement, called the causation entropy \cite{elinger2020information, elinger2021causation}, is exploited to detect the possible causal relationship between each candidate function and the time evolution of the associated state variable. The model structure is then determined by retaining those candidate functions that are demonstrated to be crucial to the underlying dynamics based on the causal inference. The causation entropy has its unique advantages in identifying model structure in the presence of indirect coupling between features and stochastic noise \cite{quinn2015directed}, which are crucial features of complex turbulent systems. In addition, with the pre-determined model structure from the causal relationship, the parameter estimation remains a quadratic optimization problem, which is much easier to solve compared with the traditional sparse identification based on a constrained optimization that involves an L1 regularization.
Second, a systematic nonlinear stochastic parameterization is built to characterize the time evolution of the unobserved state variables, aiming at capturing their statistical feedback to the observed ones \cite{chen2018conditional}. To effectively learn the details of the stochastic parameterization, which is often a computationally expensive task utilizing direct optimization algorithms, an efficient nonlinear data assimilation method is developed that exploits closed analytic formulae to sample the trajectories of the unobserved variables \cite{chen2020efficient}. These sampled trajectories are then treated as synthetic observations that allow the entire system to be fully observed, which facilitates the parameter estimation based on a simple maximum likelihood criterion.
Finally, the localization of the state variable dependence and the physics constraints with energy-conserving nonlinearity are both incorporated into the learning procedure \cite{chen2017beating, majda2012physics, chen2020learning}, which overcome the curse of dimensionality and prevent the finite time blow-up issue of the complex systems.

The rest of the paper is organized as follows. The new causality-based learning algorithm with a stochastic parameterization for complex systems with partial observations is developed in Section \ref{Sec:Method}. Examples of learning prototype complex systems are included in Section \ref{Sec:Numerics}. The paper is concluded in Section \ref{Sec:Conclusions}.

\section{The Causality-Based Data-Driven Learning Approach}\label{Sec:Method}
\subsection{Overview of the method}

Let us start with the general formulation of complex nonlinear systems \cite{vallis2017atmospheric, salmon1998lectures, kalnay2003atmospheric, majda2016introduction},
\begin{equation}\label{eq:abs_formu}
\frac{\d\mathbf{Z}}{\d t}= \boldsymbol{\Phi}(\mathbf{Z}(t)) + \boldsymbol\sigma\dot{\mathbf{W}}(t),
\end{equation}
where $\boldsymbol\Phi(\mathbf{Z}(t))$ consists of any linear and nonlinear functions of the state variable $\mathbf{Z}\in\mathbb{R}^N$, $\boldsymbol\sigma\in\mathbb{R}^{N\times d}$ is the noise amplitude and $\dot{\mathbf{W}}(t)\in\mathbb{R}^{d\times 1}$ is a white noise. For the simplicity of presentation, $d$ is assumed to be the same as $N$ and $\boldsymbol\sigma$ is assumed to be a constant diagonal matrix, which occurs in many situations. For complex systems in geophysics and fluids, $\boldsymbol{\Phi}(\mathbf{Z}(t))$ usually contains linear dispersion and dissipation, external forcing, and energy-conserving quadratic nonlinear terms. More complicated and higher order nonlinearity can be included in $\boldsymbol\Phi(\mathbf{Z}(t))$ in other applications.

Next, denote by $\mathbf{Z}=(\mathbf{X},\mathbf{Y})^\mathtt{T}$ a decomposition of the state variables, where both $\mathbf{X}$ and $\mathbf{Y}$ are multivariate with $\mathbf{X}\in\mathbb{R}^{N_1}$, $\mathbf{Y}\in\mathbb{R}^{N_2}$ and $N_1+N_2=N$. In general, $\mathbf{X}$ stands for the large-scale or resolved variables while $\mathbf{Y}$ contains a collection of medium- to small-scale variables or unresolved components. Assume one realization of the time series generated from the true underlying system is available for $\mathbf{X}$ serving as the observations while there is no observational data for $\mathbf{Y}$. For the purpose of learning the underlying physics, a library of linear and nonlinear functions that represent different combinations of the components of $\mathbf{X}$ and $\mathbf{Y}$ is pre-developed. Given an initial guess of the model structure of the observed state variables (hereafter ``model structure''), the stochastic parameterizations of the unobserved variables (hereafter ``stochastic parameterization''), and the model parameters, the learning algorithm includes an iterative procedure that alternates between three steps until the solution converges.
\begin{enumerate}
  \item Conditioned on the observed time series of $\mathbf{X}$, apply a conditional sampling to obtain a time series of the unobserved state variables $\mathbf{Y}$. The conditional sampling of $\mathbf{Y}$ can be achieved utilizing closed analytic formulae and is computationally inexpensive.
  \item Treating the sampled trajectory of $\mathbf{Y}$ as the artificial ``observations'', compute the causality-based information transfer from each candidate function in the library to the time evolution of the given state variable. Determine the model structure and the form of the stochastic parameterization based on such a causal inference.
  \item Utilize a simple maximum likelihood estimation to compute the coefficients of the above selected functions.
\end{enumerate}

For high-dimensional systems, the localization of the state variable dependence is incorporated into the causal inference such that the information transfer from only a small number of the candidate functions needs to be computed, which can mitigate the curse of dimensionality. On the other hand, physics constraints with energy conserving nonlinearity are added to the parameter estimation step, which allows the resulting model to capture the fundamental behavior of complex turbulent systems and prevents finite-time blow up of the solutions.

\begin{figure}
\centering
\includegraphics[width=1\textwidth]{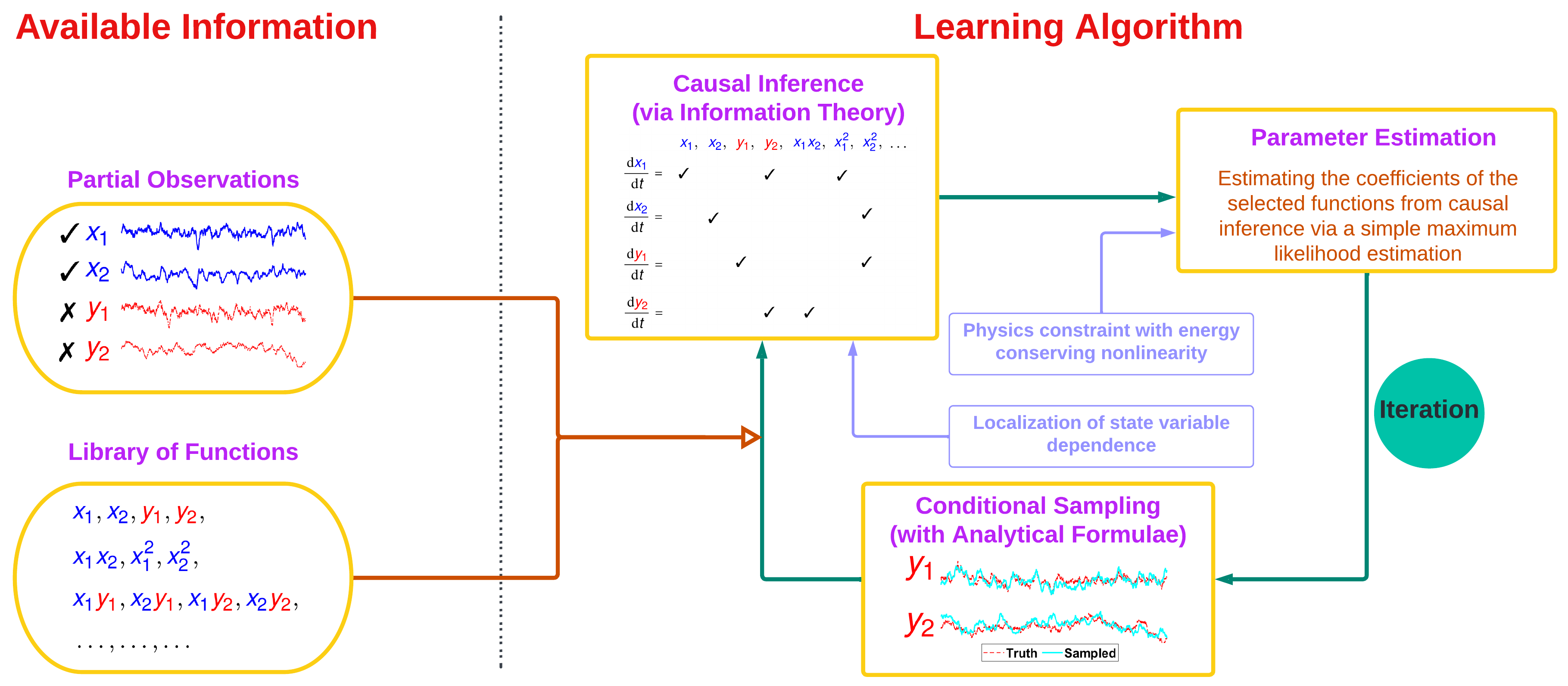}
\caption{Schematic diagram of the learning algorithm. Here, $\mathbf{X}=(x_1, x_2)$ and $\mathbf{Y}=(y_1,y_2)$.  }\label{Illustration}
\end{figure}

An overview of the proposed causality-based data-driven learning algorithm with partial observations is summarized in Figure \ref{Illustration}.

\subsection{Stochastic parameterization and conditional sampling of the unobserved state variables}
The partial observations lead to one of the fundamental challenges in efficiently learning the underlying dynamics, as the uncertainty due to the lack of observations impedes the use of simple closed formulae for the identification of the model structure as well as the estimation of the model parameters. It is worthwhile to highlight that learning the exact underlying dynamics of these unobserved variables is intrinsically very challenging if not completely infeasible for most of the complex turbulent systems, since the random noise and chaotic behavior of the signal will largely affect the observability of the system. Therefore, it is natural to build stochastic parameterizations for characterizing the unobserved variables that nevertheless can provide crucial feedback to the observed variables. To facilitate the computational efficiency of the learning process, the following stochastic parameterization structure $\mathbf{Y}$ is incorporated to the general nonlinear process of the observed state variable $\mathbf{X}$ \cite{liptser2013statistics, chen2016filtering, chen2018conditional},\begin{subequations}\label{CGNS}
\begin{align}
  \frac{\d\mathbf{X}}{\d t} &= \Big[\mathbf{A}_\mathbf{0}(\mathbf{X},t) + \mathbf{A}_\mathbf{1}(\mathbf{X},t) \mathbf{Y}(t)\Big]  + \mathbf{B}_\mathbf{1}(\mathbf{X},t)\dot{\mathbf{W}}_\mathbf{1}(t),\label{CGNS_X}\\
  \frac{\d\mathbf{Y}}{\d t} &= \Big[\mathbf{a}_\mathbf{0}(\mathbf{X},t) + \mathbf{a}_\mathbf{1}(\mathbf{X},t) \mathbf{Y}(t)\Big]  + \mathbf{b}_\mathbf{2}(\mathbf{X},t)\dot{\mathbf{W}}_\mathbf{2}(t),\label{CGNS_Y}
\end{align}
\end{subequations}
where $\mathbf{A_0}, \mathbf{a_0}, \mathbf{A_1}, \mathbf{a_1}, \mathbf{B_1}$ and $\mathbf{b_2}$ are vectors or matrices that depend nonlinearly on the observed state variable $\mathbf{X}$ and time $t$  while $\dot{\mathbf{W}}_\mathbf{1}$ and $\dot{\mathbf{W}}_\mathbf{2}$ are independent white noise sources that can have different dimensions for $\mathbf{X}$ and $\mathbf{Y}$. One key feature of the parameterization of the unobserved variable $\mathbf{Y}$ in \eqref{CGNS} is that its governing equation \eqref{CGNS_Y} is overall highly nonlinear and can produce strongly non-Gaussian statistics, but the process is conditionally linear with respect to $\mathbf{Y}$ once $\mathbf{X}$ is given. Such a family of stochastic parameterization is widely used in geophysics, climate, atmosphere, and ocean science, such as the stochastic superparameterization, dynamical super resolution, and various stochastic forecast models in data assimilation \cite{gershgorin2010test, grooms2014stochastic, branicki2018accuracy, chen2022conditional, chen2018predicting}. Since $\mathbf{Y}$ often denotes the fast, small and subgrid scale components of the system, the terms corresponding to the nonlinear self-interaction of $\mathbf{Y}$ mostly involve high frequencies and rapid fluctuations \cite{majda2003systematic}. Thus, these terms can often be effectively characterized by either simple stochastic noise \cite{majda2001mathematical, plant2008stochastic, mana2014toward, berner2017stochastic} or suitable approximations that are nonlinear functions of $\mathbf{X}$ and conditionally linear functions of $\mathbf{Y}$ \cite{mou2022efficient}. The resulting stochastic parameterization in \eqref{CGNS_Y} can successfully capture the dominant dynamics and provide similar statistics feedback to $\mathbf{X}$ as the true system. Another justification of the parameterization in \eqref{CGNS} is that many complex nonlinear systems already fit into this coupled modeling framework \cite{chen2018conditional, chen2016filtering}, including many physics-constrained nonlinear stochastic models (e.g., the noisy versions of Lorenz models,  Charney-DeVore flows, and the paradigm model for topographic mean flow interactions), a large number of stochastically coupled reaction-diffusion models in neuroscience and ecology (e.g., the FitzHugh-Nagumo models and the SIR epidemic models), and a wide class of multiscale models in turbulence and geophysical flows (e.g., the spectrum representations of the Boussinesq equations and the rotating shallow water equation). Note that the feedback from $\mathbf{Y}$ to $\mathbf{X}$ can either be through an additive function or a multiplicative one with the prefactor being an arbitrary nonlinear function of $\mathbf{X}$.

One desirable feature of the system \eqref{CGNS} is that its mathematical structure facilitates an efficient conditional sampling of the trajectory of $\mathbf{Y}$ via a closed analytic formula, which avoids sampling errors from using the particle methods and greatly accelerates the calculation.
\begin{prop}[Conditional Sampling]
For the nonlinear system \eqref{CGNS}, conditioned on one realization of the observed variable $\mathbf{X}(s)$ for $s\in[0,T]$, the optimal strategy of sampling  the trajectories associated with the unobserved variable $\mathbf{Y}$ within the same time interval satisfies the following explicit formula~\cite{chen2020efficient},
\begin{equation}\label{Sampling_Main}
  \frac{\overleftarrow{\d \mathbf{Y}}}{\d t} = \frac{\overleftarrow{\d \boldsymbol\mu_\mathbf{s}}}{\d t} - \big(\mathbf{a}_\mathbf{1} + (\mathbf{b_2}\mathbf{b}_{\mathbf{2}}^\mathtt{T})\mathbf{R}_\mathbf{f}^{-1}\big)(\mathbf{Y} - \boldsymbol\mu_\mathbf{s}) + \mathbf{b_2}\dot{\mathbf{W}}_{\mathbf{Y}}(t),
\end{equation}
where $\dot{\mathbf{W}}_{\mathbf{Y}}(t)$ is a Gaussian random noise that is independent from $\dot{\mathbf{W}}_{\mathbf{2}}(t)$ in~\eqref{CGNS}. The variables $\mathbf{R}_\mathbf{f}$ and $\boldsymbol\mu_\mathbf{s}$ are the filtering covariance and smoother mean, where the filtering and smoothing of $\mathbf{Y}$ are defined as the following conditional distributions
\begin{equation}\label{CGNS_PDF}
\begin{split}
    p(\mathbf{Y}(t)|\mathbf{X}(s),s\leq t) &\sim \mathcal{N}(\boldsymbol\mu_{\mathbf{f}}(t),\mathbf{R}_{\mathbf{f}}(t)),\\
    p(\mathbf{Y}(t)|\mathbf{X}(s), s\in[0,T])&\sim\mathcal{N}(\boldsymbol\mu_\mathbf{s}(t),\mathbf{R}_\mathbf{s}(t)),
\end{split}
\end{equation}
and the associated evolution equations are given explicitly by
\begin{subequations}\label{CGNS_Stat}
\begin{align}
  \frac{\d \boldsymbol{\mu}_{\mathbf{f}}}{\d t} &= (\mathbf{a}_\mathbf{0} + \mathbf{a}_\mathbf{1} \boldsymbol{\mu}_{\mathbf{f}}) + (\mathbf{R}_{\mathbf{f}}\mathbf{A}_{\mathbf{2}}^\mathtt{T} ) (\mathbf{B_1}\mathbf{B}_{\mathbf{1}}^\mathtt{T})^{-1} \left(\frac{\d\mathbf{X}}{\d t} - (\mathbf{A}_\mathbf{0} + \mathbf{A}_\mathbf{1}\boldsymbol{\mu}_{\mathbf{f}})\right),\label{CGNS_Stat_Mean}\\
  \frac{\d\mathbf{R}_{\mathbf{f}}}{\d t} &= \mathbf{a}_\mathbf{1} \mathbf{R}_{\mathbf{f}} + \mathbf{R}_{\mathbf{f}}\mathbf{a}_\mathbf{1}^\mathtt{T} + \mathbf{b}_{\mathbf{2}}\mathbf{b}_{\mathbf{2}}^\mathtt{T} - ( \mathbf{R}_{\mathbf{f}}\mathbf{A}_\mathbf{1}^\mathtt{T})(\mathbf{B_1}\mathbf{B}_{\mathbf{1}}^\mathtt{T})^{-1}(\mathbf{A}_\mathbf{1}\mathbf{R}_{\mathbf{f}}),\label{CGNS_Stat_Cov}\\
  \frac{\overleftarrow{\d \boldsymbol{\mu}_\mathbf{s}}}{\d t} &=  -\mathbf{a}_\mathbf{0} - \mathbf{a}_\mathbf{1}\boldsymbol{\mu}_\mathbf{s}  + (\mathbf{b_2}\mathbf{b}_{\mathbf{2}}^\mathtt{T})\mathbf{R}_{\mathbf{f}}^{-1}(\boldsymbol\mu_{\mathbf{f}} - \boldsymbol{\mu}_\mathbf{s}),\label{Smoother_Main_mu}\\
  \frac{\overleftarrow{\d \mathbf{R}_\mathbf{s}}}{\d t} &= - (\mathbf{a}_\mathbf{1} + (\mathbf{b}_{\mathbf{2}}\mathbf{b}_{\mathbf{2}}^\mathtt{T}) \mathbf{R}_{\mathbf{f}}^{-1})\mathbf{R}_\mathbf{s} - \mathbf{R}_\mathbf{s}(\mathbf{a}_\mathbf{1}^\mathtt{T} + (\mathbf{b}_{\mathbf{2}}\mathbf{b}_{\mathbf{2}}^\mathtt{T})\mathbf{R}_{\mathbf{f}})  + \mathbf{b}_{\mathbf{2}}\mathbf{b}_{\mathbf{2}}^\mathtt{T}.\label{Smoother_Main_R}
\end{align}
\end{subequations}
The notation $\overleftarrow{\d \cdot}/\d t$ in \eqref{Smoother_Main_mu}--\eqref{Smoother_Main_R} corresponds to the negative of the usual derivative, which means that both the equations are solved backward over $[0,T]$ with $(\boldsymbol\mu_\mathbf{s}(T), \mathbf{R}_\mathbf{s}(T)) = (\boldsymbol\mu_{\mathbf{f}}(T), \mathbf{R}_{\mathbf{f}}(T))$ after \eqref{CGNS_Stat_Mean}--\eqref{CGNS_Stat_Cov} have been solved forward over $[0,T]$. The starting value of the nonlinear smoother $(\boldsymbol\mu_\mathbf{s}(T), \mathbf{R}_\mathbf{s}(T))$ is the same as the end point value of the filter estimate $(\boldsymbol\mu_{\mathbf{f}}(T), \mathbf{R}_{\mathbf{f}}(T))$.
\end{prop}
Therefore, the conditional sampling formula in \eqref{Sampling_Main} allows us to recover the time series of $\mathbf{Y}$, which together with the observed trajectory of $\mathbf{X}$ forms a complete set of the time series for the entire system.

\subsection{Causal inference for discovering the model dynamics via information theory}
Given the observational time series of $\mathbf{X}$ and the sampled trajectories of $\mathbf{Y}$ from the previous step, the next task is to determine the functions in the library that are crucial to the time evolution of each state variable. Recall the collection of the state variables $\mathbf{Z}=(\mathbf{X},\mathbf{Y})^\mathtt{T}$. Denote the components of $\mathbf{Z}$ by $\mathbf{Z}=(z_1,z_2, \ldots, z_N)^\mathtt{T}$. Further denote by $\mathbf{F}=(f_1, f_2,\ldots, f_M)^\mathtt{T}$ the candidate functions in the pre-determined library. As a starting point, assume all these functions are possible candidates for the dynamics of each $z_i$, $i=1,\ldots,N$. Then, after applying a forward Euler temporal discretization scheme, the deterministic part of the starting model of $\mathbf{Z}$ has the following form:
\begin{equation}\label{Discrete_Z}
\begin{split}
  \left[\begin{array}{c}
z_{1}(t+\Delta{t}) \\
z_{2}(t+\Delta{t}) \\
\vdots \\
z_{N}(t+\Delta{t})
\end{array}\right]&=\left[\begin{array}{ccc}
\xi_{1,1} & \cdots & \xi_{1, M} \\
\xi_{2,1} & \cdots & \xi_{2, M} \\
\vdots & \ddots & \vdots \\
\xi_{N,1} & \cdots & \xi_{N, M}
\end{array}\right]\left[\begin{array}{c}
f_{1}\left(z_{1}(t), \ldots, z_{N}(t), t\right) \\
f_{2}\left(z_{1}(t), \ldots, z_{N}(t), t\right) \\
\vdots \\
f_{M}\left(z_{1}(t), \ldots, z_{N}(t), t\right)
\end{array}\right]\\
&=\boldsymbol\Xi \times \mathbf{F}\left(\mathbf{Z}(t), t\right),
\end{split}
\end{equation}
where $\boldsymbol\Xi$ is the coefficient matrix to be estimated. In general, the size of the matrix $\boldsymbol\Xi$ is quite large since $M$ is often a big number. Therefore, physics-informed sparse identification is essential to force most of the entries to be zero.

To incorporate certain physical evidence into this identification process, the following causal inference is utilized  \cite{elinger2020information}. Denote by $C_{f_{m} \rightarrow z_{n} \mid\left[\mathbf{F} \backslash {f}_{m}\right]}$ the causation entropy of $f_{m}(t)$ on $z_{n}(t+\Delta{t})$ conditioned on all $\mathbf{F}$ except $f_m$, which allows to explore the composition of $z_{n}(t+\Delta{t})$ that comes solely from $f_{m}(t)$. If such a causation entropy is zero (or practically nearly zero), then $f_{m}(t)$ does not contribute any information to  $z_{n}(t+\Delta{t})$ and the associated parameter $\xi_{n,m}$ is set to be zero. By computing such a causation entropy for different $m=1,\ldots, M$ and $n=1,\ldots,N$, a sparse causation entropy matrix is reached which indicates if each entry of $\boldsymbol\Xi$ should be estimated. Note that $f_m$ is supposed to be a non-constant function, and the constant terms are always assumed to exist. Then a simple maximum likelihood estimation based on a  quadratic optimization can be easily applied to determine the actual values of those nonzero entries in $\boldsymbol\Xi$.

The causation entropy $C_{f_{m} \rightarrow z_{n} \mid\left[\mathbf{F} \backslash {f}_{m}\right]}$ is defined as follows,
\begin{equation}\label{Causation_Entropy}
\begin{split}
  C_{f_{m} \rightarrow z_{n} \mid\left[\mathbf{F} \backslash {f}_{m}\right]} &= H(z_{n}|\left[\mathbf{F} \backslash {f}_{m}\right]) - H(z_n|\left[\mathbf{F} \backslash {f}_{m}\right], f_{m})\\
  &=H(z_{n}|\left[\mathbf{F} \backslash {f}_{m}\right]) - H(z_n|\mathbf{F}).
\end{split}
\end{equation}
In \eqref{Causation_Entropy}, $H(\cdot|\cdot)$ is the conditional entropy, which is related to the Shannon's entropy $H(\cdot)$ and the joint entropy $H(\cdot,
\cdot)$. They are defined as:
\begin{equation*}
\begin{split}
  H(X) &= -\int_x p(x)\log(p(x))\d x,\\
  H(Y| X) &= -\int_x\int_y p(x,y)\log(p(y|x))\d y\d x,\\
  H(X,Y) &= -\int_x\int_y p(x,y)\log(p(x,y))\d y\d x,
\end{split}
\end{equation*}
where $p$ is the associated PDF that is typically determined by a histogram from the given time series (assuming the ergodicity).
The difference on the right-hand side of \eqref{Causation_Entropy} naturally represents the contribution of $f_m$ to $z_n$.

Note that other information measurements have also been applied for the causality detection, such as the transfer entropy \cite{barnett2009granger, vicente2011transfer} and the directed information \cite{massey1990causality, kramer1998directed}. These methods in general work well. Nevertheless, the causation entropy in \eqref{Causation_Entropy} has its unique advantages in identifying model structure in the presence of indirect coupling between features and stochastic noise \cite{quinn2015directed}, which are crucial features of complex turbulent systems.

\medskip

\noindent \textbf{The calculation of the causation entropy.}\\
The direct calculation of the causation entropy in \eqref{Causation_Entropy} is nontrivial and is computationally expensive. In fact, reconstructing the exact PDFs from a given time series is a very challenging issue. The kernel density estimation (KDE) \cite{schreiber2000measuring}, the box-counting algorithm \cite{paluvs2001synchronization} and many other direct estimation methods suffer from the curse of dimensionality. Some alternative methods have been proposed, such as the k-nearest
neighbors \cite{kraskov2004estimating, murphy2012machine}, which can mitigate the issue to some extent but may remain to be complicated. Since determining the model structure only depends on if the causation entropy is zero or not, rather than its exact value, the following properties will facilitate the calculation of the causation entropy in high dimensions.
\begin{prop}[Chain rule]
The conditional entropy can be represented by the Shannon's entropy and the joint entropy via the following chain rule:
\begin{equation}
H(Y | X) = H(X, Y)-H(X).
\end{equation}
\end{prop}
\begin{prop}[Gaussian approximation]
If $p \sim \mathcal{N}(\boldsymbol{\mu}, \mathbf{R})$ satisfies a $s$-dimensional Gaussian distribution, then the Shannon's entropy has the following explicit form
\begin{equation}
H(p)=\frac{s}{2}\left(1+\ln (2 \pi)\right)+\frac{1}{2} \ln\left( \operatorname{det}(\mathbf{R})\right),
\end{equation}
where `det' is the matrix determinant.
\end{prop}
With these properties in hand, the practical calculation of the causation entropy can be the following.
\begin{prop}[Practical calculation of the causation entropy]
By approximating all the joint and marginal distributions as Gaussians, the  causation entropy can be computed in the following way:
\begin{equation}
\label{Entropy Guassians}
\begin{split}
C_{Z \rightarrow X | Y} &=H(X | Y)-H(X | Y, Z) \\
& = H(X,Y) - H(Y) - H(X,Y,Z) + H(Y,Z)\\
& = \frac{1}{2} \ln(\operatorname{det}(\mathbf{R}_{XY}))-\frac{1}{2} \ln(\operatorname{det}(\mathbf{R}_{Y})) \\&\qquad\qquad\qquad- \frac{1}{2} \ln(\operatorname{det}(\mathbf{R}_{XYZ})) +\frac{1}{2} \ln(\operatorname{det}(\mathbf{R}_{YZ})),
\end{split}
\end{equation}
where $\mathbf{R}_{XYZ}$ denotes the covariance matrix of the state variables $(X,Y,Z)^\mathtt{T}$ and similar for other covariances.
\end{prop}
It is worthwhile to remark that the Gaussian approximation may lead to certain errors in computing the causation entropy if the true distribution is highly non-Gaussian. Nevertheless, the primary goal here is not to obtain the exact value of the causation entropy. Instead, it suffices to detect all the index pairs $(m,n)$ in $\boldsymbol\Xi$, associated with which the causation entropy $C_{f_{m} \rightarrow z_{n} \mid\left[\mathbf{F} \backslash {f}_{m}\right]}$ is nonzero (or practically above a small threshold value). In most of the applications, if a significant causal relationship is detected in the higher order moments, then very likely it exists in the Gaussian approximation as well. This allows us to efficiently determine the structure of the sparse matrix of $\boldsymbol\Xi$, where the exact values of the nonzero entries will be calculated via a simple maximum likelihood estimation.

Similar to \eqref{eq:abs_formu}, the system by retaining only the functions corresponding to the nonzero causation entropy entries can be written as
\begin{equation}\label{eq:abs_formu_identified}
\frac{\d\mathbf{Z}}{\d t}= \widetilde{\boldsymbol{\Phi}}(\mathbf{Z}(t)) + \boldsymbol\sigma\dot{\mathbf{W}}(t).
\end{equation}
Further denote by $\boldsymbol\Theta$ the collection of the parameters to be estimated, which correspond to the nonzero entries in $\boldsymbol\Xi$.

\subsection{Parameter estimation via a simple maximum likelihood estimation }
Consider a temporal discretization of \eqref{eq:abs_formu_identified} using the Euler-Maruyama scheme \cite{gardiner2004handbook},
\begin{equation}\label{discrete_eqn}
  \mathbf{Z}^{j+1} = \mathbf{Z}^j + \widetilde{\boldsymbol\Phi}(\mathbf{Z}^j) \Delta{t} + \boldsymbol\sigma\boldsymbol\varepsilon^j\sqrt{\Delta{t}},
\end{equation}
where $j$ is the index in time, $\Delta{t}$ is a small time step, and $\boldsymbol\varepsilon^j$ is an independent and identically distributed (i.i.d.) standard multidimensional Gaussian random number. Denote by $\mathbf{z}^{j}$ the given numerical value of $\mathbf{Z}^{j}$ from observations. Further denote by $\mathbf{M}^j \boldsymbol\Theta + \mathbf{s}^j := \mathbf{z}^j + \widetilde{\boldsymbol\Phi}(\mathbf{z}^j) \Delta{t}$, namely the deterministic part on the right hand side of \eqref{discrete_eqn} evaluated at $\mathbf{z}^j$. Here the nonlinear functions are included in $\mathbf{M}^j$ and $\mathbf{s}^j$, where the former appears as the multiplicative prefactor of the parameters $\boldsymbol\Theta$  while the latter appears on its own such as the first term --- $\mathbf{Z}^j$ --- on the right-hand side of \eqref{discrete_eqn}. Due to the Euler-Maruyama approximation, the one-step time evolution from $\mathbf{Z}^j$ to $\mathbf{Z}^{j+1}$ is approximated by a linear function within such a short time interval. Therefore, the likelihood can be computed based on a Gaussian distribution,
\begin{equation}
  \mathcal{N} (\boldsymbol\mu^j, \boldsymbol\Sigma) = C |\boldsymbol\Sigma|^{-\frac{1}{2}} \exp \left(- \frac{1}{2} (\mathbf{z}^{j+1} - \boldsymbol\mu^j)^\mathtt{T} (\boldsymbol\Sigma)^{-1} (\boldsymbol\mu^{j+1} - \boldsymbol\mu^j)\right),
\end{equation}
where the mean and the covariance are given by $\boldsymbol\mu^j = \mathbf{M}^j \boldsymbol\Theta + \mathbf{s}^j$ and $\boldsymbol\Sigma = \boldsymbol\sigma\boldsymbol\sigma^\mathtt{T}\Delta{t}$, respectively. Note that $\boldsymbol\sigma$ has been assumed to be a diagonal constant matrix, and therefore $\boldsymbol\Sigma$ does not depend on $j$. Taking a logarithm operation to cancel the exponential function and summing up the likelihood over the entire time period yield
\begin{equation}\label{eq: SI_obj}
  \mathcal{L} =  \frac{1}{2} \sum_{j=1}^J  (\mathbf{z}^{j+1} - \mathbf{M}^j \boldsymbol\Theta - \mathbf{s}^j)^\mathtt{T}  (\boldsymbol\Sigma)^{-1} (\mathbf{z}^{j+1} - \mathbf{M}^j \boldsymbol\Sigma - \mathbf{s}^j)  - \frac{J}{2} \log |\boldsymbol\Sigma| ,
\end{equation}
where $J+1=\lfloor T/\Delta{t}\rfloor$ with $\lfloor\cdot\rfloor$ being the floor function that rounds down the result to the nearest integer.
To find the minimum of~$\mathcal{L}$, it is sufficient to find the zeros of $\frac{\partial\mathcal{L}}{\partial \boldsymbol\Theta} =0$ and $\frac{\partial \mathcal{L}}{\partial \boldsymbol\Sigma} =0$, which leads to
\begin{subequations}\label{eq: SI_Equation_R_Theta}
    \begin{align}
        \boldsymbol\Sigma &= \frac{1}{J} \sum_{j=1}^J (\mathbf{z}^{j+1} - \mathbf{M}^j \boldsymbol\Theta - \mathbf{s}^j)(\mathbf{z}^{j+1} - \mathbf{M}^j\boldsymbol\Theta - \mathbf{s}^j)^\mathtt{T} , \label{eq: SI_Equation_R_Theta1}\\
        \boldsymbol\Theta &= \mathbf{D}^{-1} \mathbf{c},\label{eq: SI_Equation_R_Theta2}
    \end{align}
\end{subequations}
where
\begin{equation}\label{eq: SI_aux_physics_constraint}
	\mathbf{D} = \sum_{j=1}^J (\mathbf{M}^j)^\mathtt{T} \boldsymbol{\Sigma}^{-1}\mathbf{M}^j \quad \text{and} \quad \mathbf{c} = \sum_{j=1}^J(\mathbf{M}^j)^\mathtt{T} \boldsymbol{\Sigma}^{-1}(\mathbf{z}^{j+1} - \mathbf{s}^j).
\end{equation}
The equations in \eqref{eq: SI_Equation_R_Theta} are solved by first setting $\boldsymbol\Theta=\mathbf{0}$ in finding $\boldsymbol\Sigma$ in \eqref{eq: SI_Equation_R_Theta1} via essentially the quadratic variation, and then plugging in the result into \eqref{eq: SI_Equation_R_Theta2} and \eqref{eq: SI_aux_physics_constraint} to obtain $\boldsymbol\Theta$.

\subsection{Physics constraints}
Physics constraints, meaning the conservation of energy in the quadratic nonlinear terms, are important properties in many complex turbulent systems and appear in most of the classical geophysics and fluid models \cite{majda2012physics, harlim2014ensemble}. The physics constraints prevent finite-time blow up of the solutions and facilitate a skillful medium- to long-range forecast. Therefore, taking into account the physics constraints and other constraints is crucial for the learning algorithm, especially in the parameter estimation step. These constraints, together with other constraints,  can in general be represented in the following way:
\begin{equation}
	\mathbf{H} \boldsymbol\Theta = \mathbf{g},
\end{equation}
where $\mathbf{H}$ and $\mathbf{g}$ are constant matrices. To incorporate these constraints, the Lagrangian multiplier method is applied, which modifies the  objective function in \eqref{eq: SI_obj},
\begin{equation}\label{eq: objective_constraint}
\begin{gathered}
	\mathcal{L} = \frac{1}{2} \sum_{j=1}^J (\mathbf{z}^{j+1} - \mathbf{M}^j \boldsymbol\Theta - \mathbf{s}^j)^\mathtt{T}  (\boldsymbol\Sigma)^{-1} (\mathbf{z}^{j+1} - \mathbf{M}^j \boldsymbol\Theta - \mathbf{s}^j)\\ - \frac{J}{2}\log |\boldsymbol\Sigma^{-1}|  + \boldsymbol{\lambda}^\mathtt{T}  (\mathbf{H} \boldsymbol\Theta - \mathbf{g}),
\end{gathered}
\end{equation}
where $\boldsymbol\lambda$ is the Lagrangian multiplier.
The solution to the minimization problem with the new objective function \eqref{eq: objective_constraint} is given as follows,
\begin{subequations}\label{eq: SI_Equation_R_Theta_lambda}
    \begin{align}
        \boldsymbol\Sigma &= \frac{1}{J} \sum_{j=1}^J  (\mathbf{z}^{j+1} - \mathbf{M}^j \boldsymbol\Theta - \mathbf{s}^j)(\mathbf{z}^{j+1} - \mathbf{M}^j\boldsymbol\Theta - \mathbf{s}^j)^\mathtt{T} \label{eq: SI_Equation_R_Theta_lambda1}\\
        \boldsymbol\lambda &= \left(\mathbf{H} \mathbf{D}^{-1} \mathbf{H}^\mathtt{T} \right)^{-1} (\mathbf{H} \mathbf{D}^{-1}\mathbf{c} - \mathbf{g}), \\
        \boldsymbol\Theta &= \mathbf{D}^{-1} \left( \mathbf{c} - \mathbf{H}^\mathtt{T} \boldsymbol\lambda \right),
    \end{align}
\end{subequations}
where $\mathbf{D}$ and $\mathbf{c}$ are defined in~\eqref{eq: SI_aux_physics_constraint}.

\subsection{Localization of the state variable dependence}
Localization of the state variable dependence is a typical feature in many complex turbulent systems for both the modeling of large-scale dynamics and the stochastic parameterizations. On one hand, the high-dimensional stochastic ordinary differential equations (SDEs) are usually obtained as a result of the spatial discretization of a stochastic partial differential equation. The advection, diffusion, and dispersion are all local operators, which implies that each state variable in the SDEs interacts with only the nearby few states \cite{majda2003introduction, vallis2017atmospheric}. On the other hand, the stochastic parameterizations of the states at the subgrid scales also depend only on the nearby corresponding large-scale state variables \cite{grabowski2004improved, gagne2020machine, chattopadhyay2020data}. Besides, the idea of localization is widely utilized in data assimilation and prediction \cite{bergemann2010localization, anderson2007exploring, janjic2011domain}.

The localization of the state variables is incorporated into the proposed learning algorithm at both the causal detection and the conditional sampling steps. The necessity of the localization in the causal detection is to mitigate the curse of dimensionality. In fact, when the dimension of the system becomes large, the number of functions in the library that includes different combinations of the state variables increases exponentially. As a result, 
the cost of computing the causation entropy for all these functions also shoots up. The localization, which requires to compute the causation entropy of only those functions that involve the local interactions of the state variables, can overcome the curse of dimensionality. Next, the stochastic parameterizations of each component of the subgrid variable $\mathbf{Y}$ in \eqref{CGNS_Y}, denoted by $y_{i,j}$, depends only on the associated large-scale observed state variable of $\mathbf{X}$, namely $x_i$ and $x_{i\pm1},\ldots x_{i\pm s}$ with $s$ being a small positive integer. This leads to a block covariance matrix of $\mathbf{R}_\mathbf{f}$ in \eqref{Sampling_Main} when carrying out the conditional sampling to recover the trajectory of $\mathbf{Y}$. In other words, the giant covariance matrix not only becomes a sparse one but can be divided into several low-dimensional blocks, which are then solved in a parallel way. This significantly facilitates the learning algorithm in applying to high-dimensional systems, which are otherwise difficult to handle due to the heavy computational burden in storing and solving the full covariance matrix.

\section{Test Examples}\label{Sec:Numerics}
In this section, three nonlinear chaotic or turbulent systems are utilized for testing the learning algorithm developed in the previous section. The first test model is a three-dimensional low-order system, which is mainly used as a proof-of-concept and to display the detailed procedure of the method. The other two models are spatially-extended multiscale systems, which are adopted to understand the skill of the model identification and the stochastic parameterizations.

In all these experiments, physics constraints are incorporated into the learning algorithm. Localization is adopted in the second and third experiments, which facilitates the reduction of the computational cost. In all the experiments, the total length of the observation is $500$ units while the numerical integration time step is $\Delta{t}=0.001$ such that there are in total $500,000$ points in each observational time series.
\subsection{Lorenz 1984 model}
The first test example is a low-dimensional chaotic system, known as the Lorenz 1984 (L-84) model, which is a simple analog of the global atmospheric circulation \cite{vallis2017atmospheric, salmon1998lectures}. It has the following form \cite{lorenz1984formulation, lorenz1984irregularity}:
\begin{equation}\label{Lorenz84}
\begin{split}
  \frac{\d x}{\d t} &=  {-(y^2+z^2)} - a( x-f)  + \sigma_x\dot{W}_x,\\
  \frac{\d y}{\d t} &=  -bxz +  x y - y + g  + \sigma_y\dot{W}_y,\\
  \frac{\d z}{\d t} &=  b x y +  x z- z + \sigma_z\dot{W}_z.
\end{split}
\end{equation}
In \eqref{Lorenz84}, the zonal flow $x$ represents the intensity of the mid-latitude westerly wind current, and a wave component exists with $y$ and $z$ representing the cosine and sine phases of a chain of vortices superimposed on the zonal flow. Relative to the zonal flow, the wave variables are scaled so that
$x^2 + y^2 + z^2$ is the total scaled energy. These equations can be derived as a Galerkin truncation of the two-layer quasigeostrophic potential vorticity equations in a channel. The additional stochastic noise represents the interactions between these resolved scale variables and the unresolved ones.

The following parameters are utilized for the test here:
\begin{equation}\label{L84_Coefficients}
  a = \frac{1}{4}, \quad b = 4, \quad f = 8,\quad g = 1, \quad \mbox{and}\quad \sigma_x = \sigma_y = \sigma_z = 0.1,
\end{equation}
which are the standard parameters that create chaotic behavior \cite{lorenz1984formulation}. In fact, $f > 1$ is a necessary condition for the zonal flow becoming unstable, forming steadily progressing vortices, while $g > 0$ triggers the chaotic behavior of the entire system.

\subsubsection{The experiment setup}
In this experiment, $y$ and $z$ are taken as the observed variables while the observation of $x$ is not directly available. Note that the L-84 model \eqref{Lorenz84} automatically fits into the framework \eqref{CGNS} with $\mathbf{X}=(y,z)^\mathtt{T}$ and $\mathbf{Y}=x$. This fact, together with the small system noises $\sigma_x = \sigma_y = \sigma_z = 0.1$ utilized here, allows the learning algorithm to have a potential to fully recover the system, including the unobserved process, as the contribution from the deterministic part of the dynamics is only weakly polluted by the system noises.

It is natural to incorporate all the linear and quadratic nonlinear functions of $y$ and $z$, namely $y$, $z$, $y^2$, $z^2$, and $yz$, into the library of the candidate functions. This mimics the general form of the geophysical flows, the nonlinearity of which is dominated by the quadratic terms. In addition, the linear and conditional linear functions of $x$, namely $x$, $xy$, $xz$, as well as the constant forcing term are included in the library. Essentially, all the quadratic nonlinear functions of the three state variables, except $x^2$ that breaks the structure of \eqref{CGNS}, are contained in the library of the candidate functions. To further increase the complexity of the library, the cubic terms that satisfy \eqref{CGNS} are also added, which contain the quadratic terms of $y$ or $z$ multiplying $x$ but not the quadratic or cubic functions of $x$ itself.

A random and complicated initial model structure is utilized to start the iterative algorithm,
\begin{equation}\label{Lorenz84_InitialGuess}
\begin{split}
  \frac{\d x}{\d t} &=  y^2 - z^2 + 2 + (y^2 - z^2)x  + \sigma_x\dot{W}_x,\\
  \frac{\d y}{\d t} &=  - y - 2y^2 + z^2 + 1 + (- y - 8z - yz)x  + \sigma_y\dot{W}_y,\\
  \frac{\d z}{\d t} &=  - z + z^2 - yz + (8y + z + z^2)x + \sigma_z\dot{W}_z.
\end{split}
\end{equation}
The initial values of the noise coefficients in the observed processes, namely $\sigma_y$ and $\sigma_z$, are chosen to be $1$. In fact, the initial values of these two parameters won't affect the learning algorithm, as they will converge to the truth within one iteration step based on the quadratic variation \eqref{eq: SI_Equation_R_Theta_lambda1}.
On the other hand,  $\sigma_x$ is not uniquely determined from the algorithm since the effect due to the increase of $\sigma_x$ can be completely offset by decreasing the coefficients in front of $x$ in the observed processes. Therefore, if the primary goal is to learn the dynamics of the observed variables with a reasonable parameterization of the unobserved ones, then an arbitrary value of $\sigma_x$ can be used. For the simplicity of the study here, $\sigma_x=0.1$ is set to be known. It is also worthwhile to remark that, as the quadratic variation of the unobserved variable in general cannot be directly updated by the conditional sampling, a change of variable to normalize such a diffusion coefficient is often adopted to update such a parameter if the coefficients in front of $x$ in the observed processes are known \cite{beskos2005exact}.

\subsubsection{Results}

\begin{figure}[!htb]
\centering
\includegraphics[width=1\textwidth]{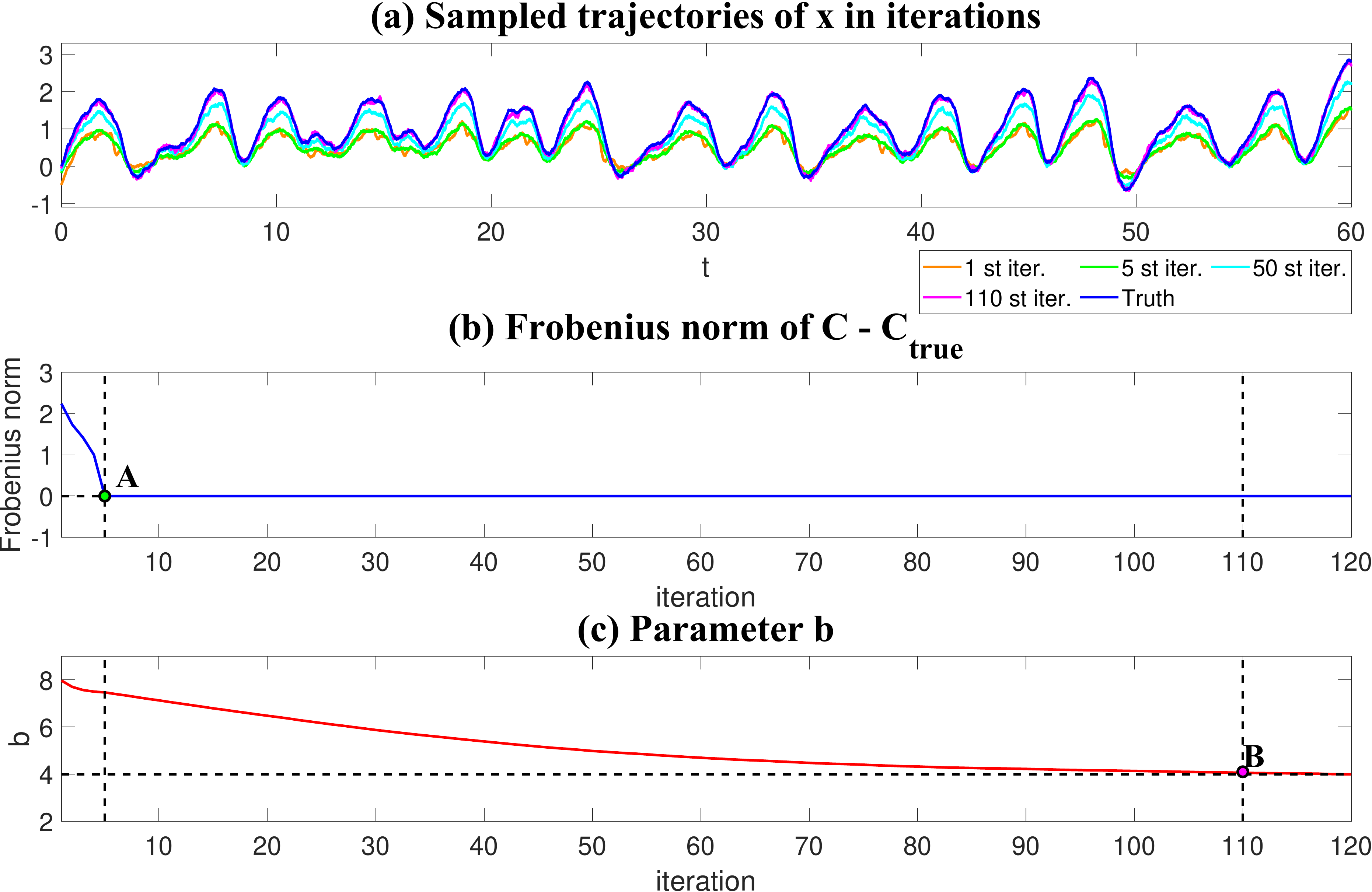}
\caption{Iterative procedure of learning the L-84 model with partial observations $(y,z)^\mathtt{T}$. Panel (a): The sampled trajectory of the unobserved variable $x$ at the 1st, the 5th, the 50th, and the 110th iterations, where there are in total $120$ iterations. The random initial guess of the model structure is shown in \eqref{Lorenz84_InitialGuess}. Panel (b): The Frobenius norm of $\mathbf{C} - \mathbf{C}_{\text{true}}$ as a function of iterations, where $\mathbf{C}$ is the causation entropy matrix indicator at each iteration step with entries being either $0$ or $1$ and $\mathbf{C}_{\text{true}}$ is the causation entropy matrix indicator corresponding to the true model \eqref{Lorenz84}. Panel (d): The updates of the parameter $b$. The points $\mathbf{A}$ and $\mathbf{B}$ show the iteration at 5 and 110 steps. The model structure is identified correctly after the 5th iteration step, while the parameters converge to the truth at around the 110th step.}
\label{Sampling history}
\end{figure}

Figure \eqref{Sampling history} displays the detailed procedure of the learning algorithm. Panel (a) shows the sampled trajectory of the unobserved variable $x$ at the 1st, the 5th, the 50th, and the 110th iterations. It is seen that the sampled trajectory converges to the truth as the number of iterations increases, indicating that the learning process eventually recovers the unobserved trajectory and identifies the model structure.

To better understand the iterative procedure, Panel (b) of Figure \eqref{Sampling history} shows the convergence of the model structure towards the truth. Here, a causation entropy matrix indicator $\mathbf{C}$ is introduced, which is of size $N\times M$, where $N=3$ is the dimension of the system while $M$ is the total number of candidate functions. The matrix $\mathbf{C}$ has the same structure as $\boldsymbol\Xi$ in \eqref{Discrete_Z} except that $\mathbf{C}$ is a logical matrix with entries being either $0$ or $1$. If the causation entropy associated with a specific term exceeds the pre-determined threshold (which is $10^{-3}$ here), then the corresponding entry in $\mathbf{C}$ is set to be $1$, meaning that the term should be maintained in the dynamics. Then the Frobenius norm of $\mathbf{C} - \mathbf{C}_{\text{true}}$ is computed, where $\mathbf{C}_{\text{true}}$ is the causation entropy matrix indicator corresponding to the true model \eqref{Lorenz84}. It is seen that despite the large gap in the initial random guess of the model structure, there are only $5$ terms (corresponding to Frobenius being $2.2361$) that are mismatched after $1$ iteration (the first point in the curve) and the correct structure is reached after merely $5$ iterations. Note that, at the $5$th iteration, the sampled trajectory of $x$ (green) in Panel (a) is still far from the truth. This is because, although the model structure is already perfectly identified, additional iterations are still required for the parameters to converge. One desirable feature observed in Figure \eqref{Sampling history} is that the model structure does not change after the $5$-th iteration, but only the parameters are updated that simultaneously provide the improved sampled trajectory of $x$. The final parameter values after $120$ iterations are shown in Table \ref{Table:L84}, which are almost indistinguishable from the truth. Here, the model \eqref{Lorenz84} is rewritten in the following form for the convenience of comparing the parameters displayed in Table \ref{Table:L84},
\begin{equation}\label{Lorenz84_Identified}
\begin{split}
  \frac{\d x}{\d t} &=  {\theta_{yy}^{x} y^2+\theta_{zz}^{x} z^2 } +\theta_x^x x + \theta_1^x + \sigma_x\dot{W}_x,\\
  \frac{\d y}{\d t} &=  \theta_{xz}^y xz + \theta_{xy}^y xy + \theta_y^y y + \theta_1^y + \sigma_y\dot{W}_y,\\
  \frac{\d z}{\d t} &=  \theta_{xy}^z x y + \theta_{xz}^z x z + \theta_z^z z + \theta_1^y + \sigma_z\dot{W}_z.
\end{split}
\end{equation}
\begin{table}[h]
\centering
\begin{tabular}{|l|c|c|c|c|c|c|}
\hline & $\theta_{x}^{x}$ & $\theta_{y}^{y}$ & $\theta_{z}^{z}$ & $\theta_{yy}^{x}$ & $\theta_{zz}^{x}$ & $\theta_{xz}^{y}$ \\
\hline \text { Truth } & -0.2500 & -1.0000 & -1.0000 & -1.0000 & -1.0000 & -4.0000 \\
\text { Identified }   & -0.2680 & -0.9987 & -1.0076 & -0.9993 & -1.0061 & -3.9956 \\
\hline  & $\theta_{xy}^{y}$ & $\theta_{xy}^{z}$ & $\theta_{xz}^{z}$ & $\theta_{1}^{x}$ & $\theta_{1}^{y}$& $\theta_{1}^{z}$\\
\hline \text { Truth } & 1.0000 & 4.0000 & 1.0000 & 2.0000 & 1.0000& 0.0000\\
\text { Identified }  & 0.9993 & 3.9956 & 1.0061 & 2.0223 & 0.9939& 0.0053 \\
\hline
\end{tabular}\caption{Comparison of the parameters in the true system \eqref{Lorenz84} and those in the identified model. Since the Frobenius norm of $\mathbf{C} - \mathbf{C}_{true}$ converges to zero, the identified model has exactly the same structure and the truth. They are both rewritten in the form of \eqref{Lorenz84_Identified} for the convenience of comparing the model parameters. }\label{Table:L84}
\end{table}

\begin{figure}[!htb]
\centering
\includegraphics[width=1\textwidth]{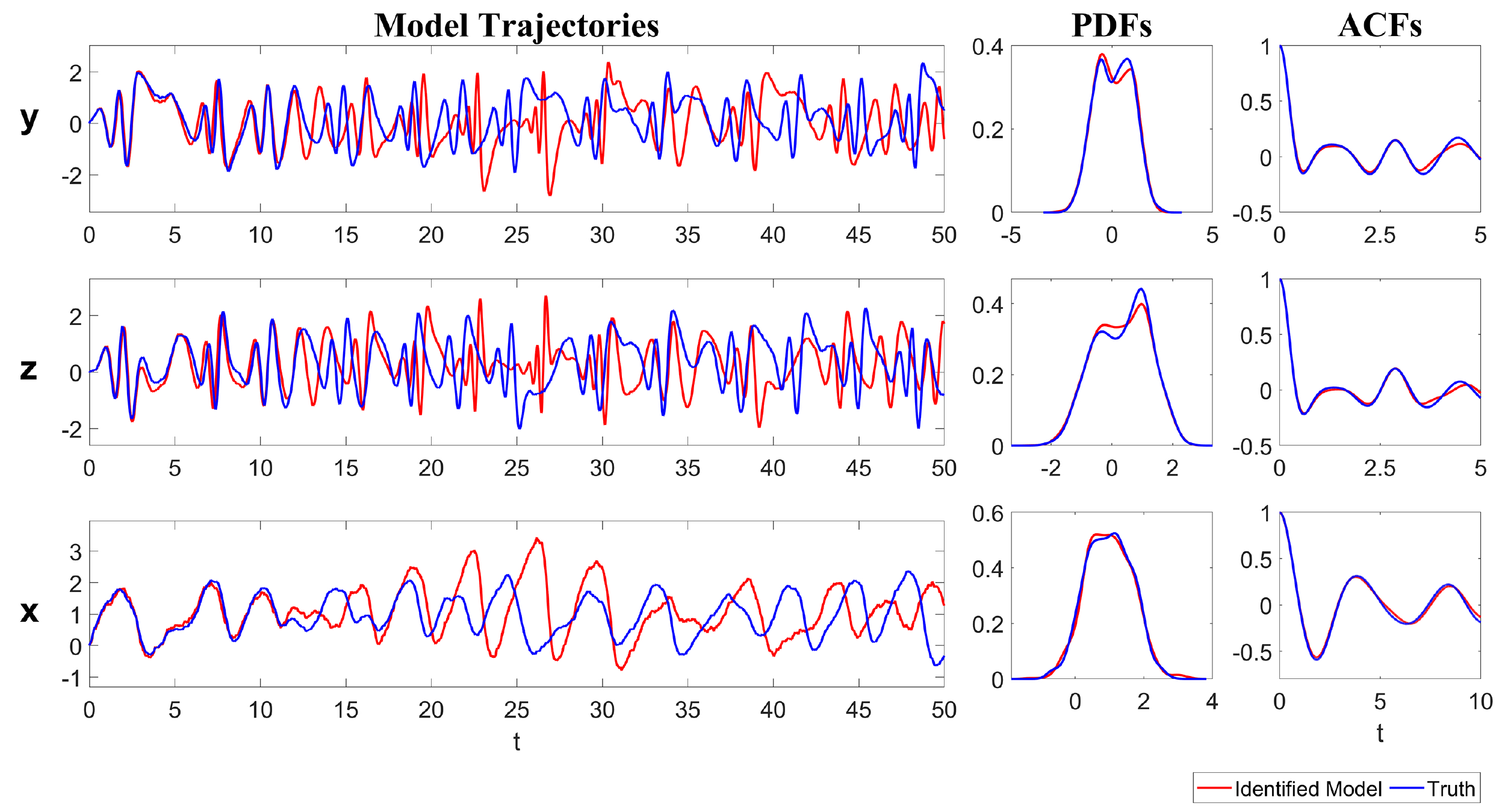}
\caption{Comparison of the model trajectories and the associated statistics using the true model (blue) and the identified model (red). The two statistics utilized are the PDF and the temporal autocorrelation function (ACF). Note that due to the chaotic nature of the system, the two curves do not expect to have a one-to-one point-wise match between each other. Instead, a qualitative similarity between the trajectories from the truth and the identified model is evidence to show the accuracy of the identified model.}
\label{Solutions}
\end{figure}

Figure \ref{Solutions} shows the model simulations and the associated statistics using the identified model together with the estimated parameters. Due to the chaotic nature of the system, the trajectories from the true and the identified models do not expect to have a one-to-one point-wise match between each other. Nevertheless, a qualitative similarity between the trajectories from the truth and the identified model is observed, which indicates the accuracy of the identified model. This is further confirmed by the nearly perfect recovery of the two statistics: the PDF and the temporal autocorrelation function (ACF). These facts conclude the skill of the algorithm based on this simple chaotic example.

\subsection{Two-layer Lorenz 1996 models}
The two-layer Lorenz 1996 (L-96) model \cite{wilks2005effects, lorenz1996predictability} is a conceptual representation of geophysical turbulence that is commonly used as a tested for data assimilation and parameterization in numerical weather forecasting \cite{majda2009mathematical,grooms2014stochastic,grooms2014stochastic2,majda2014new}. The model mimics a coarse discretization of atmospheric flow on a latitude circle. It supports complex wave-like and chaotic behavior, and the two-layer structure schematically depicts the interactions between small-scale fluctuations and large-scale motions. The stochastic version of the model subject to additive noise forcing reads
\begin{subequations}\label{L96_v1}
	\begin{align}
		\frac{\d u_{i}}{\d t} &=  - u_{i - 1}\left(u_{i - 2}-u_{i + 1}\right)-u_{i} + f - \frac{h c_i}{J} \sum^{J}_{j=1} v_{i, j}  \notag \\
&\qquad\qquad\qquad\qquad\qquad\qquad + \sigma_{u_i} \dot W_{u_i}, \quad i = 1, \dots, I, \label{L96_v1_u} \\
		\frac{\d v_{i, j}}{\d t} &=  -b c_i v_{i, j + 1}\left(v_{i, j + 2} - v_{i, j - 1}\right) - c_i v_{i, j} + \frac{h c_i}{J} u_{i}   \notag\\
&\qquad\qquad\qquad\qquad\qquad\qquad + \sigma_{v_{i, j}} \dot W_{v_{i, j}}, \quad j=1, \dots, J, \label{L96_v1_v}
\end{align}
\end{subequations}
where $I$ denotes the total number of large-scale variables and $J$ is the number of small-scale variables corresponding to each large-scale variable. In \eqref{L96_v1}, $f$, $h$, $c_i$, $b$, $\sigma_{u_i}$ and $\sigma_{v_{i,j}}$ are given scalar parameters while $\dot W_{u_i}$ and $\dot W_{v_{i, j}}$ are independent white noises. The large-scale variables $u_i$s are periodic in $i$ with $u_{i + I} = u_{i - I} = u_{i}$. The corresponding small-scale variables~$v_{i, j}$s are periodic in $i$ with $v_{i + I, j} = v_{i - I, j} = v_{i, j}$ and satisfy the following conditions in $j$: $v_{i, j + J} = v_{i + 1, j}$, and $v_{i, j - J} = v_{i - 1, j}$.

\begin{figure}[!htb]
\centering
\includegraphics[width=1\textwidth]{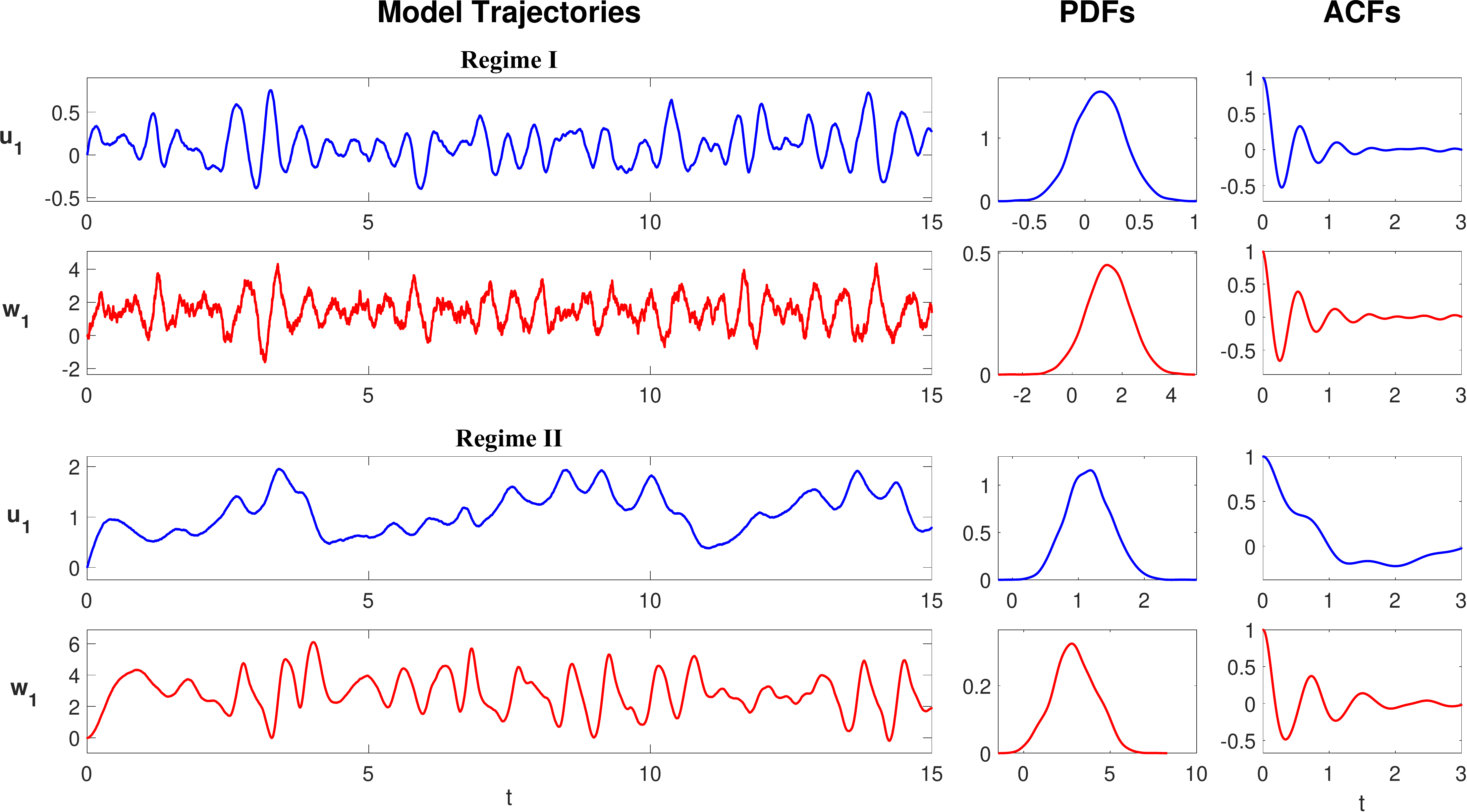}
\caption{Dynamical regimes of the Two-layer L-96 model \eqref{L96_v1}, including the model trajectories, PDFs and ACFs of variables $u_i$ and $w_i$ at $i=1$. Here $w_i$ is defined as  $w_i=\sum_{j=1}^J v_{i,j}$. }
\label{L96_parameter}
\end{figure}

Two dynamical regimes are considered here as the truth. They share most of the parameters:
\begin{equation}\label{L96_Coefficients_1}
I=20, \quad J=4, \quad c_{i}=2+0.7 \cos (2 \pi i / I), \quad b=2, \quad f=4, \quad \sigma_{u_{i}}=0.05,
\end{equation}
but  they are differed by $h$ and $\sigma_{v_{i, j}}$:
\begin{equation}\label{L96_Coefficients_2}
\begin{split}
  \mbox{Regime I:}&\qquad h = 4.0 \qquad \mbox{and}\qquad\sigma_{v_{i, j}}=1.00\\
  \mbox{Regime II:}&\qquad h = 1.5\qquad\mbox{and}\qquad \sigma_{v_{i, j}}=0.05.
\end{split}
\end{equation}
The model trajectories and statistics of these two dynamical regimes are shown in Figure \ref{L96_parameter}. For the convenience of discussing the behavior of the two layers, a new single variable $w_i=\sum_{j=1}^J v_{i,j}$ is introduced, which describes the total variabilities in the second layer. It can be seen that there is no scale separation between $u_i$ and $v_i$ in regime I since the ACFs oscillate and decay in a similar fashion. On the other hand, $u_i$ tends to occur in a slower time scale compared with $v_i$ in Regime II, leading to multiscale features. It is also worthwhile to mention that the coefficient $c_{i}$ is spatially varying, which gives an inhomogeneous spatial pattern of the system.

\subsubsection{The experiment setup}

Here the $u_i$s for $i = 1,2,\dots, I$ are the observed variables and all the $v_{i,j}$s for different $i$ and $j$ are the unobserved ones. Since only the time series of $u_i$s are provided while the structure of the true model is unavailable, the number $J$ is unknown to us. Therefore, a natural way to build a suitable model is to incorporate stochastic parameterizations into the processes of $u_i$s. Each of such a stochastic parameterization $w_i$ takes into account the total contributions of the associated $v_{i,j}$ to a specific $u_i$, which is effectively $w_i = \sum^{J}_{j=1} v_{i, j}$. Therefore, the target model has $2I$ dimensions, with $I$ state variables being the observed ones and the remaining $I$ variables representing the stochastic parameterizations.

Due to the high dimensionality of the problem, the size of the library consisting of the candidate functions will become extremely large if all possible linear and nonlinear functions up to a certain order are considered. Nevertheless, since the main components of the dynamics, such as the advection and diffusion, involve only local interactions, it is natural to consider such localizations in building the library of the candidate functions. To this end, for each $u_{i}$, only the terms involving its adjacent variables $u_{i-1},u_{i-2},u_{i+1}$ and $u_{i+2}$ are utilized to construct the candidate functions. The nonlinearity considered here is up to the quadratic terms. In addition, the contribution from the stochastic parameterization needs to be included. One of the simplest choices is to augment the library by one additive term $w_i$ and one multiplicative term $u_iw_i$, where $w_i$ itself is driven by a hidden process representing the stochastic parameterization. Note that other more complicated nonlinear interactions between the state variables $u_i$s and the stochastic parameterizations can be easily included in the library. But for the parsimony of the model, only these two related terms are utilized here.
The set of the candidate functions for $u_i$ is then given by a vector $\mathbf{F}_{u_i}$, which includes $23$ terms:
\begin{equation}\label{Candidate_u_i_L96}
\begin{gathered}
  u_{i},\ u_{i-1},\ u_{i-2},\ u_{i+1},\ u_{i+2},\ u_{i}^2,\ u_{i-1}^2,\ u_{i-2}^2,\ u_{i+1}^2,\\
   u_{i+2}^2,\ u_{i}u_{i-1},\ u_{i}u_{i-2},\ u_{i}u_{i+1},\ u_{i}u_{i+2},\ u_{i-1}u_{i-2},\ u_{i-1}u_{i+1},\\
  u_{i-1}u_{i+2},\ u_{i-2}u_{i+1},\ u_{i-2}u_{i+2},\ u_{i+1}u_{i+2},\ 1,\ w_i,\ u_iw_i.
\end{gathered}
\end{equation}
Clearly, the candidate functions in the library allow rich features to appear in the dynamics, such as the diffusion and other quadratic nonlinear interactions that were not in the true system.
On the other hand, only $4$ terms are included in the library for each $w_i$, given by another vector $\mathbf{F}_{w_i}$,
\begin{equation}\label{Candidate_v_i_L96}
  u_i,\ u_i^2,\ 1,\ w_i.
\end{equation}
This allows for a simple form of the stochastic parameterization. Nevertheless, the nonlinear terms $u_i^2$ in  $\mathbf{F}_{w_i}$ and $u_iw_i$ in  $\mathbf{F}_{u_i}$ satisfy the physics constraints.

The initial guess of the model is constructed as follows:
\begin{subequations}\label{L96_v1_guess}
	\begin{align}
		\frac{\d u_{i}}{\d t} &=  - u_{i - 1}\left(u_{i - 2}-u_{i + 1}\right) + \left(u_{i+1}^2-u_{i-1}u_{i}\right)  + \left(u_{i}u_{i+1}- u_{i-1}^2\right) -u_{i} + f  \notag \\
&\qquad\qquad\qquad\qquad\qquad\qquad\qquad - \frac{h c_i}{J} w_i + \sigma_{u_i} \dot W_{u_i}, \quad i = 1, \dots, I, \label{L96_v1_u_guess} \\
		\frac{\d w_i}{\d t} &= h c_{i} u_{i} -  w_{i} + \sigma_{w_{i}} \dot W_{w_{i}},\quad i = 1, \dots, I, \label{L96_v1_v_guess}
\end{align}
\end{subequations}
where the terms $u_{i+1}^2-u_{i-1}u_{i}$ and $u_{i}u_{i+1}- u_{i-1}^2$ are two pairs of local quadratic advection added beyond those in the true system. Both the pairs of local quadratic advection satisfy the physical constraints.

\begin{figure}
\centering
\includegraphics[width=1\textwidth]{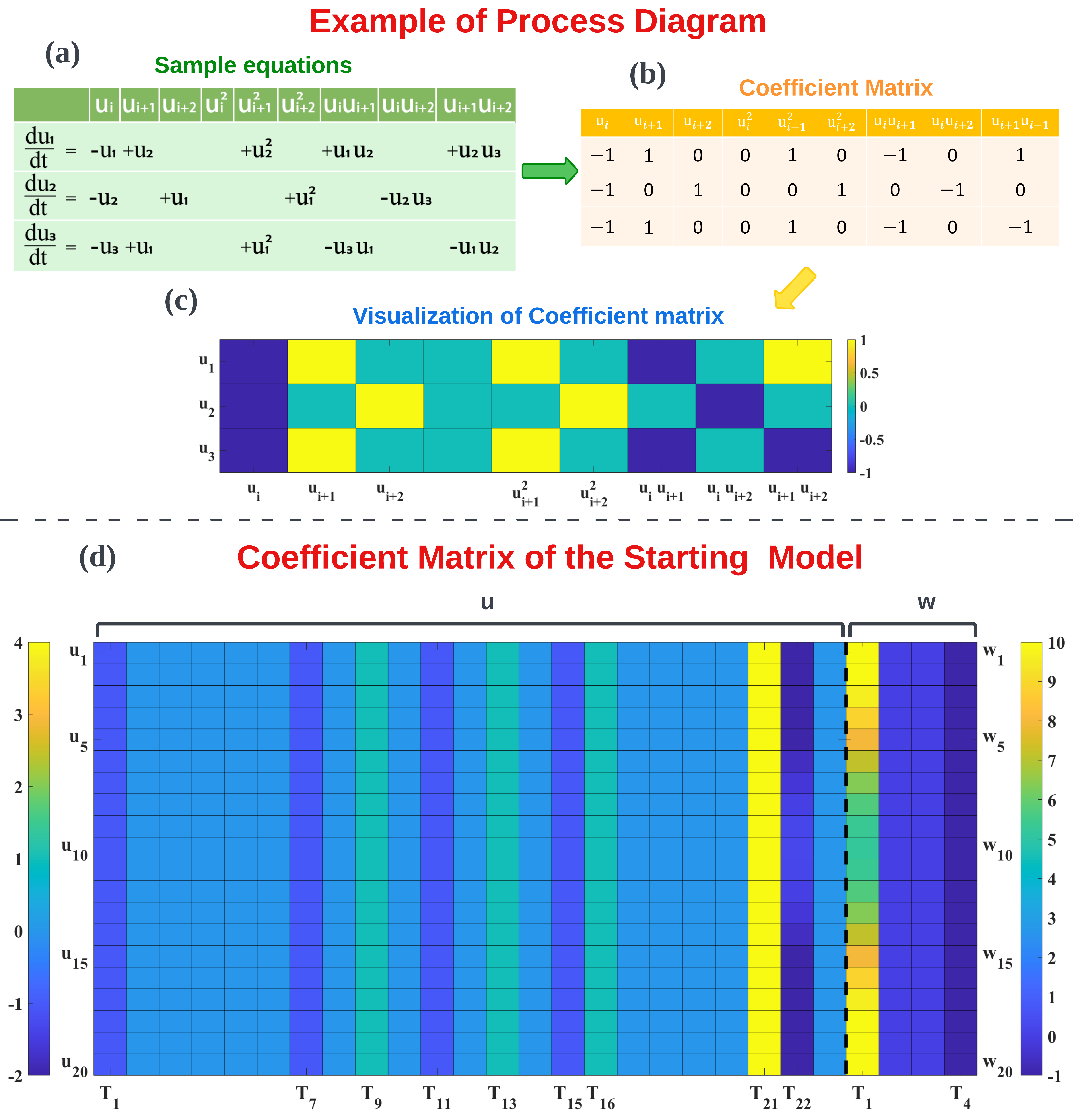}
\caption{Using a visualization diagram to represent the identified model structure and parameters. Panels (a)--(c): a low-dimensional ODE system as a simple example to illustrate the visualization of the coefficient matrix. Panel (d): the coefficient matrix of the two-layer L-96 model corresponding to the starting model \eqref{L96_v1_guess}. In the big coefficient matrix, the $i$-th row represent the right hand side of the equations of $u_i$ and $w_i$, where $u_i$ and $w_i$ contain $24$ and $4$ terms, respectively. The order of these terms in the figure is the same as that in \eqref{Candidate_u_i_L96} and \eqref{Candidate_v_i_L96}. The parameter values are indicated by the colors. The $\mathbf{T}_j$ under the $u$ and $w$ parts stands for the $j$-th term in \eqref{Candidate_u_i_L96} and \eqref{Candidate_v_i_L96}, namely $\mathbf{F}_{u_i}$ and $\mathbf{F}_{w_i}$, respectively.}
\label{coefficient_matrix}
\end{figure}

In the following, a visualization diagram is utilized to represent the identified model structure and parameters. Figure \ref{coefficient_matrix} includes an illustration of the visualization diagram for the starting model. Panels (a)--(c) show the general representation of the sparse coefficient matrix. Panel (d) corresponds to the starting model \eqref{L96_v1_guess}. In the big coefficient matrix, the $i$-th row represent the right hand side of the equations of $u_i$ and $w_i$, where $u_i$ and $w_i$ contain $23$ and $4$ terms, respectively. The order of these terms in the figure is the same as that in \eqref{Candidate_u_i_L96} and \eqref{Candidate_v_i_L96}.  The parameter values are indicated by the colors.

\subsubsection{Results}
\begin{figure}[h]
\centering
\includegraphics[width=1\textwidth]{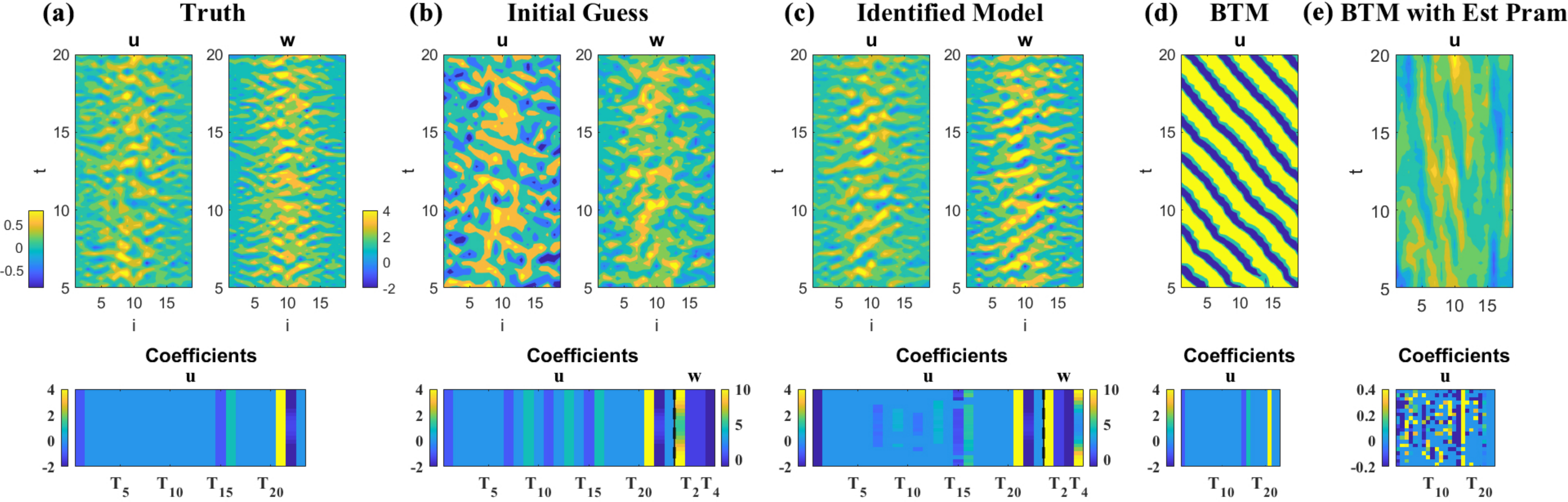}
\caption{Identifying the two-layer L-96 model in Regime I. Different columns show the truth, the initial guess, the identified model, and the bare truncation models (BTMs). The first row displays the spatiotemporal pattern of both $u_i$ and $w_i$. Note that different random number seeds are utilized in different columns, and there is no point-wise correspondence between different patterns. The focus is only on the overall structure. The second row shows the coefficient matrix in each model, as was described in Figure \ref{coefficient_matrix}. The threshold value in determining the causation matrix for the identified model with stochastic parameterization is $0.001$ and that for the BTM model (Panel (e)) is $0.0001$ since the same threshold as the former leads to an even worse result.}
\label{L96_comparison}
\end{figure}

\begin{figure}[h]
\centering
\includegraphics[width=1\textwidth]{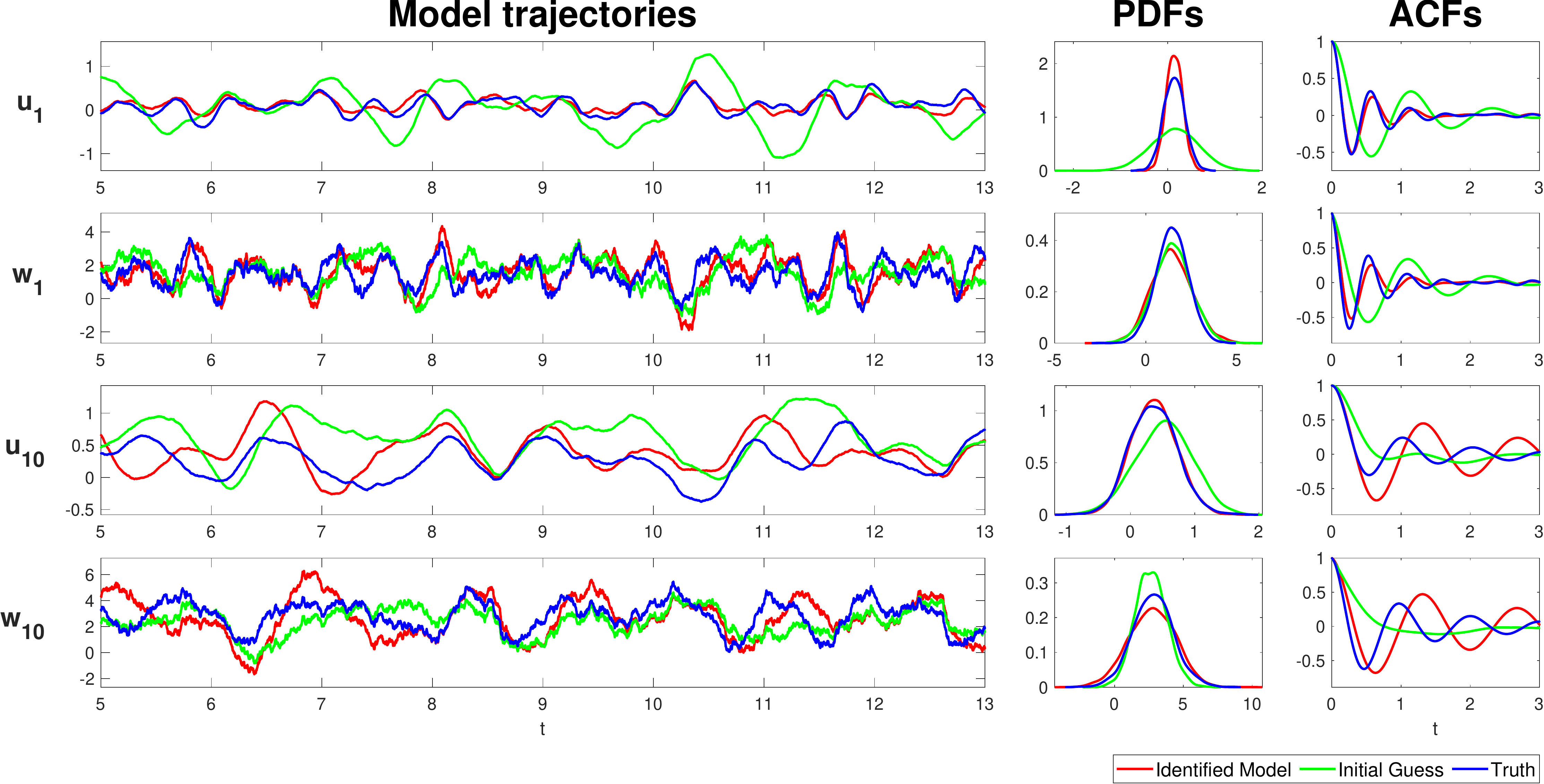}
\caption{Comparison of the model trajectories, the ACFs and the PDFs of the truth, the initial guess of the model, and the identified model in Regime I, at $i=1$ and $i=10$. Similar to Figure \ref{L96_comparison}, different random number seeds are utilized in different models, and there is no point-wise correspondence between the trajectories from different models. Only the qualitative features are the useful information obtained from these trajectories. The threshold value in determining the causation matrix is $0.001$.}
\label{L96_comparison2}
\end{figure}

Let us start with Regime I, in which $u_i$ and $v_i$ lie on the same time scale. A threshold of $0.001$ is utilized in determining if each entry in the causation entropy matrix should be retained.
After $50$ iterations, the results converge, where the identified model is shown in Panel (c) of Figure \ref{L96_comparison}. It is seen that the identified model is qualitatively similar to the truth (Panel (a)). In particular, the spatial inhomogeneous structure in $u_i$ is recovered. In addition, despite the chaotic behavior, the weakly eastward propagation of the individual waves and the westward propagation of the wave envelope in the spatiotemporal pattern of $u_i$ are both captured by the identified model. Figure \ref{L96_comparison2} compares the model trajectories and statistics. It is clear that the trajectories of $u_i$ generated from the identified model (which uses a different random number generator from the truth) are qualitatively similar to the truth, and the statistics are much more accurate than the initial guess.

Next, it is important to understand the role of the stochastic parameterization. To this end, the so-called bare truncation model (BTM) is adopted for the purpose of comparison, which only retains the dynamics of $u_i$ but completely omits the equations of $w_i$. That is, only \eqref{L96_v1_u} is utilized, where all $v_{i,j}$s are set to be zero. Thus, the BTM has a dimension of $I$. It is shown in Panel (d) of Figure \ref{L96_comparison} that if the same parameters as in the true system are adopted for the BTM, then the wave patterns become much more regular than the truth due to the lack of perturbations from the small scales. Even by incorporating a parameter estimation to the BTM, the spatiotemporal pattern of the BTM is different from the truth (Panel (e)). This indicates the important role of $w_i$s in the original system, especially in such a case that there is no clear scale separation in the true system. Therefore, incorporating a stochastic parameterization is essential to characterize its effect.
One important finding is that the stochastic parameterization $w_i$ actually recovers the combined contribution of all the $v_{i,j}$ for $j = 1,\ldots, J$, according to Figures \ref{L96_comparison}--\ref{L96_comparison2}. This is good evidence that shows the role of the one-dimensional stochastic parameterization in replacing the $J$-dimensional small-scale features in the true model. Note that the coefficients of $w_i u_i$ are zero in the identified model, which implies that the single additive term $w_i$ is sufficient to parameterize the total contribution of the small-scale feedbacks.

\begin{figure}[!htb]
\centering
\includegraphics[width=1\textwidth]{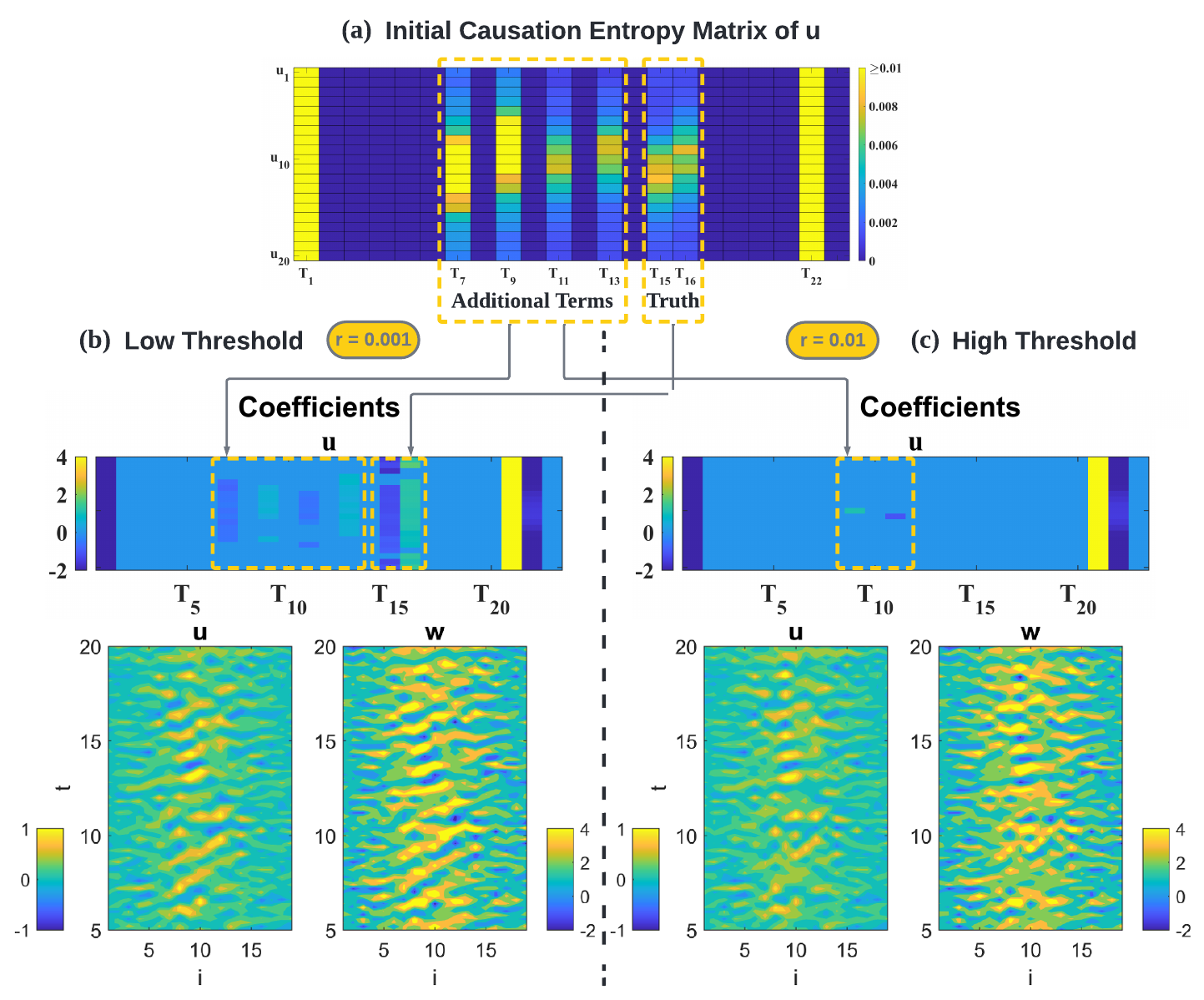}
\caption{Comparison of the identified model of the two-layer L96 model in Regime I using different threshold values for determining the causation entropy matrix. Panel (a): the causation entropy matrix of $u_i$s corresponding to the initially guessed structure except for column $\mathbf{T_{21}}$ representing trivial constant terms. Panel (b): the coefficient matrix for $u_i$ based on the low threshold $r= 0.001$ and the corresponding spatial-temporal patterns. Panel (c): The coefficient matrix for $u$ based on the high threshold $r=0.01$ and the corresponding spatial-temporal patterns.}
\label{L96_thrd}
\end{figure}

Yet, this is perhaps one undesirable feature in the identified model. That is, there are several additional terms remained in the identified model (e.g., $\mathbf{T_7}$, $\mathbf{T_9}$, $\mathbf{T_{11}}$ and $\mathbf{T_{13}}$), which do not appear in the perfect system. This is because of the specific threshold used here to determine if each candidate function should be kept. The threshold value is $r=0.001$, which is relatively low. Therefore, it is natural to repeat the learning process but increase the threshold. To this end, a higher threshold $0.01$ is utilized, and the results are shown in Figure \ref{L96_thrd}. It is seen that not only those additional terms but also the advection terms disappear with this high threshold. The reason that the advection terms rather than the damping $u_i$, the forcing $f$, and the feedback from $w_i$ are chosen to be eliminated by the learning algorithm is because of its relatively weak role in the original dynamics in this special regime. In fact, according to Panel (c) of Figure \ref{L96_thrd}, the spatiotemporal pattern of $u_i$ remains similar to the truth by a quick glance. In particular, the spatial inhomogeneity, the strengths of the signal and the frequency at each fixed spatial grid point all resemble the truth. Yet, by a careful comparison with the truth, the weakly eastward propagation of the wave no longer exists, which is obviously due to the ignorance of the advection. Therefore, the comparison here indicates that the threshold value is a useful criterion in determining the importance of different terms in the identified system.

\begin{figure}[!htb]
\centering
\includegraphics[width=1\textwidth]{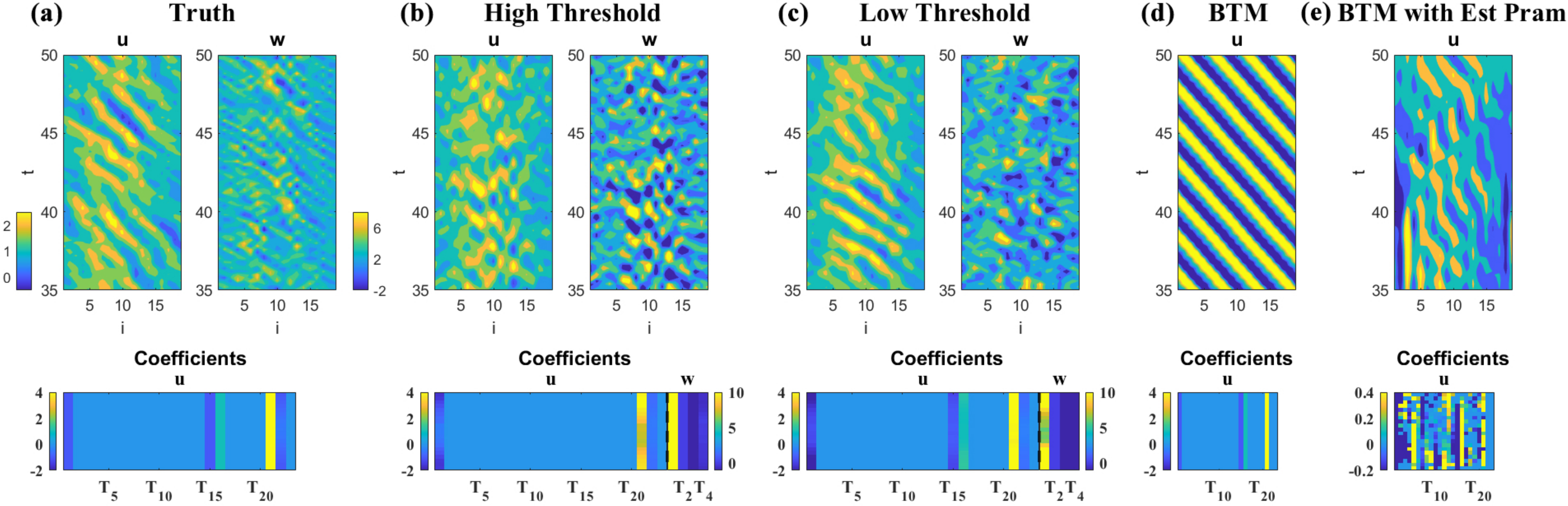}
\caption{Similar to Figure \ref{L96_comparison} but for Regime II. Note that instead of repeating the column for the initial guess, the column for the high threshold case $r=0.01$ is shown instead.}
\label{L96_threshold}
\end{figure}

\begin{figure}[!htb]
\centering
\includegraphics[width=1\textwidth]{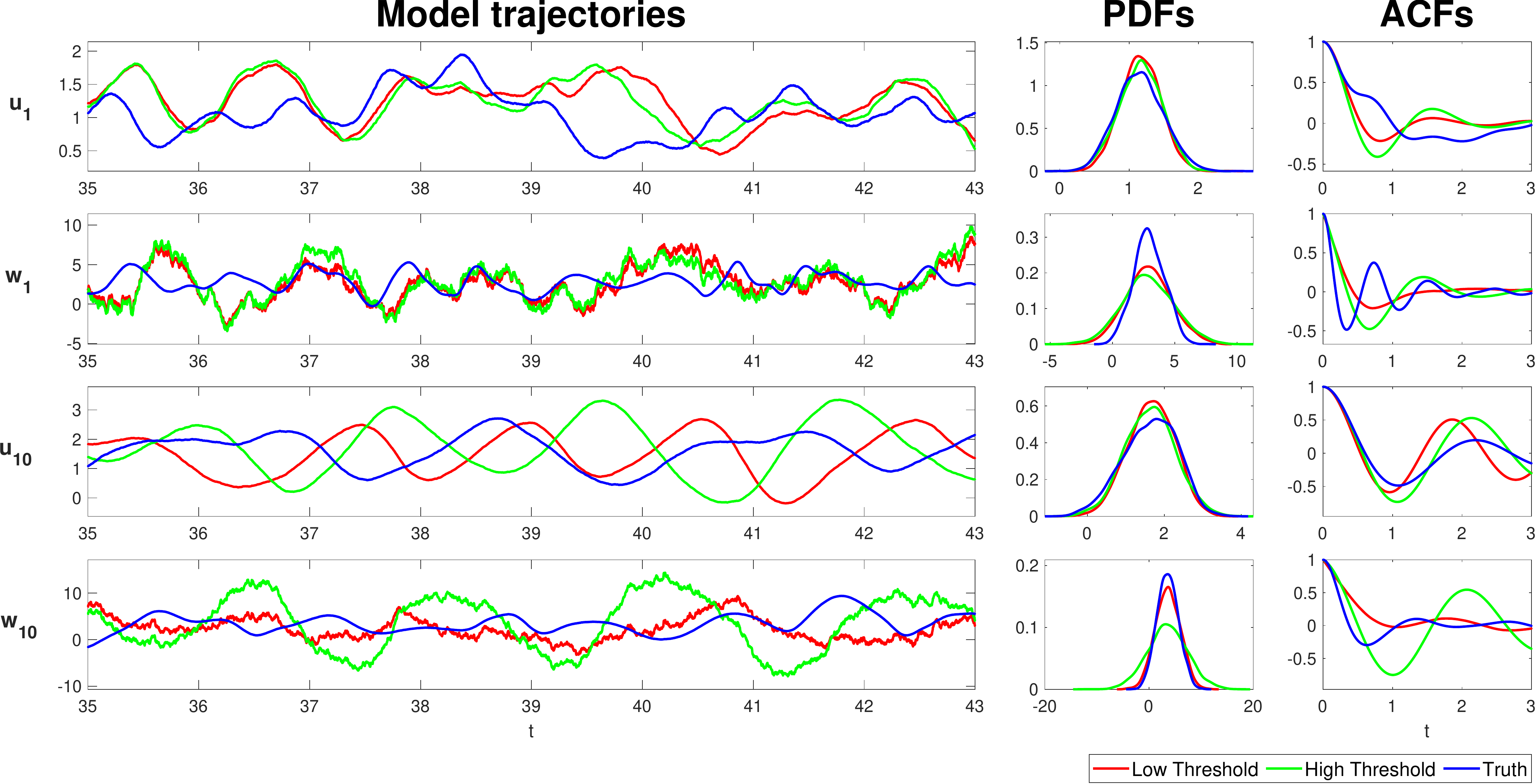}
\caption{Similar to Figure \ref{L96_comparison2} but for Regime II. Note that instead of repeating the curves for the initial guess, the curves for the high threshold case $r=0.01$ are shown instead. }
\label{L96_threshold2}
\end{figure}

Finally, Figure \ref{L96_threshold} shows the truth and the identified model in Regime II. Different from Regime I, where $u_i$ and $w_i$ lie in the same time scale and the signal of $u_i$ is quite chaotic, a clearer wave propagation pattern is observed in the spatiotemporal pattern in Regime II. This indicates the more significant role played by the advection, as the feedback from the small-scale variables $w_i$ becomes less dominant. The identified model with the threshold being $r=0.001$ again reproduces most features of the underlying dynamics. On the other hand, if a higher threshold $r=0.01$ is utilized, then the advection is again eliminated by the learning algorithm. However, in this dynamical regime, the spatiotemporal pattern of $u_i$ without the advection becomes quite distinct from the truth due to its missing the wave propagations. Figure \ref{L96_threshold2} shows the model trajectories, the ACFs, and the PDFs, which confirm the skill of the identified model in recovering the dynamical and statistical features of the observed variables $u_i$s. Although there are certain gaps between the truth and the stochastic parameterized process $w_i$, all that is important for $w_i$ in the identified model is its feedback to $u_i$ but not its own exact dynamics.

\subsection{A stochastically coupled FitzHugh–Nagumo (FHN) model}
The last test model for the learning algorithm is the following stochastically coupled FitzHugh-Nagumo (FHN) model. The FHN model is a prototype of an excitable system, which describes the activation and deactivation dynamics of a spiking neuron \cite{lindner2004effects}. Stochastic versions of the FHN model with the notion of noise-induced limit cycles were widely studied and applied in the context of stochastic resonance \cite{treutlein1985noise, lindner2000coherence, longtin1993stochastic, wiesenfeld1994stochastic} and its spatially extended version has also attracted much attention as a noisy excitable medium \cite{neiman1999noise, hempel1999noise, hu2000phase, jung1998noise}. By exploiting a finite difference discretization to the diffusion term, the stochastically coupled FHN model in the lattice form is given by
\begin{equation}\label{FHN}
\begin{split}
  \epsilon \frac{\d u_i}{\d t} &= \left(d_u(u_{i+1}+u_{i-1}-2u_{i}) + u_i - \frac{1}{3}u_i^3 - v_i\right) +\sqrt{\epsilon}\delta_1\dot{W}_{u_i},\\
   \frac{\d v_i}{\d t} &= \big(u_i + a\big)  + \delta_2\dot{W}_{v_i},\qquad i = 1,\ldots, N,
\end{split}
\end{equation}
where the parameter $a > 1$ is required in order to guarantee that the system has a global attractor in the absence of noise and diffusion. The random noise is able to drive the system above the threshold level of global stability and triggers limit cycles intermittently. The time scale ratio $\epsilon$ is much smaller than one, implying that the $u_i$s are the fast and the $v_i$s are the slow variables. The model in \eqref{FHN} is equipped with the spatial periodic boundary conditions.

The following parameters are utilized to generate the observational time series:
\begin{equation}\label{FHN_Coefficients}
  \epsilon = 0.01,\qquad \delta_1 = 0.2,\qquad \delta_2 = 0.1,\qquad d_u = 10, \qquad \mbox{and} \qquad a = 1.05.
\end{equation}
The number of the spatial grid points is $N=40$.
With these parameters, the structure of the solution exhibits strong spatial coherent structures in $u_i$. See Panel (a) of Figure \ref{FHN_comparison}.

\subsubsection{The experiment setup}

The time series of $u_i$ for $i = 1, 2, 3, \dots, N$ are observed while there is no direct observation for $v_{i}$.
Similar to the setup in the L-96 model, define a vector $\mathbf{F}_{u_i}$ that contains $28$ candidates functions for each $u_i$:
\begin{equation}\label{FHN_Candidate_u}
\begin{gathered}
  u_{i},\ u_{i-1},\ u_{i-2},\ u_{i+1},\ u_{i+2},\ u_{i}^2,\ u_{i-1}^2,\ u_{i-2}^2,\ u_{i+1}^2,\ u_{i+2}^2,\,u_{i}u_{i-1}, \\
u_{i}u_{i-2},\ u_{i}u_{i+1},\ u_{i}u_{i+2},\ u_{i-1}u_{i-2},\ u_{i-1}u_{i+1},\ u_{i-1}u_{i+2},\ u_{i-2}u_{i+1},\\
 u_{i-2}u_{i+2},\ u_{i+1}u_{i+2},\ u_i^3,\ u_{i-1}^3,\ u_{i-2}^3,\ u_{i+1}^3,\ u_{i+2}^3,\ 1,\ v_i,\ u_i v_i,
\end{gathered}
\end{equation}
 and another vector $\mathbf{F}_{v_i}$ that includes $4$ candidate functions for each $v_i$:
\begin{equation}\label{FHN_Candidate_v}
 u_i,\ u_i^2,\ 1,\ v_i.
\end{equation}
When these candidate functions are defined, localization has been utilized as the process of each $u_i$ only depends on those terms within the nearby grid points. Nevertheless, compared with the true system, many more additional advection, diffusion and other quadratic nonlinear terms appear in the candidate functions, which are not in the stochastically coupled FHN system \eqref{FHN}. The study here contains the following two experiments.

In the first experiment, the same model structure is taken in the initial guess as in the truth \eqref{FHN}. However, very different parameters $d_u = 0.5$ and $\delta_2 = 0.4$ are utilized in the initial guess, which lead to a completely distinct spatiotemporal structure compared with the truth. In fact, as is shown in Panel (b) of Figure \ref{FHN_comparison}, the coherent structure of the model corresponding to the initial guess only appears intermittently and happens in local regions. In addition, the time series of $u_i$ at a fixed spatial location $i$ has a much higher frequency than that of the truth.

In the second experiment, two kinds of the additional terms are added to the starting model: (i) the local quadratic advection satisfying the physical constraints such as $u_{i+1}^2$ and $u_{i}u_{i-1}$, $u_{i-1}^2$ and $u_{i}u_{i+1}$, and (ii) $-v_i u_i$ in the equations of $u_i$ and the quadratic term of $u_i^2$ in the equations of $v_i$. In addition, $d_u = 0.5$ and $\delta_2 = 0.1$ are adopted for the initial guess of the model. The initial model reads:
\begin{subequations}\label{FHN_guess}
	\begin{align}
		\epsilon \frac{\d u_{i}}{\d t} &=  \frac{1}{2}\left(u_{i+1}+u_{i-1}\right) -\frac{1}{3}u_i^2 - v_i + \left(u_{i+1}^2-u_{i-1}u_{i}\right)+ \left(u_{i}u_{i+1}- u_{i-1}^2\right) \notag \\
		&\quad  + \left(u_{i}u_{i+2}- u_{i-2}^2\right)  + \epsilon u_i v_i + \sqrt{\epsilon}\delta_1\dot{W}_{u_i}, \quad i = 1, \dots, I, \label{FHN_u_guess} \\
		\frac{\d v_i}{\d t} &= \big(u_i - u_i^2 + a\big)  + \delta_2\dot{W}_{v_i},\qquad i = 1,\ldots, N,\label{FHN_v_guess}
\end{align}
\end{subequations}
where the terms $\left(u_{i+1}^2-u_{i-1}u_{i}\right)$, $\left(u_{i}u_{i+1}- u_{i-1}^2\right)$ and $\left(u_{i}u_{i+2}- u_{i-2}^2\right)$ are the three pairs of local quadratic advection added in the starting model satisfying the physics constraints. Panel (b) of Figure \ref{FHN_advection} shows the spatiotemporal patterns of such a starting model. It has much faster temporal frequencies than the truth despite the coherent structures.


\subsubsection{Results}
\begin{figure}[!htb]
\centering
\includegraphics[width=1\textwidth]{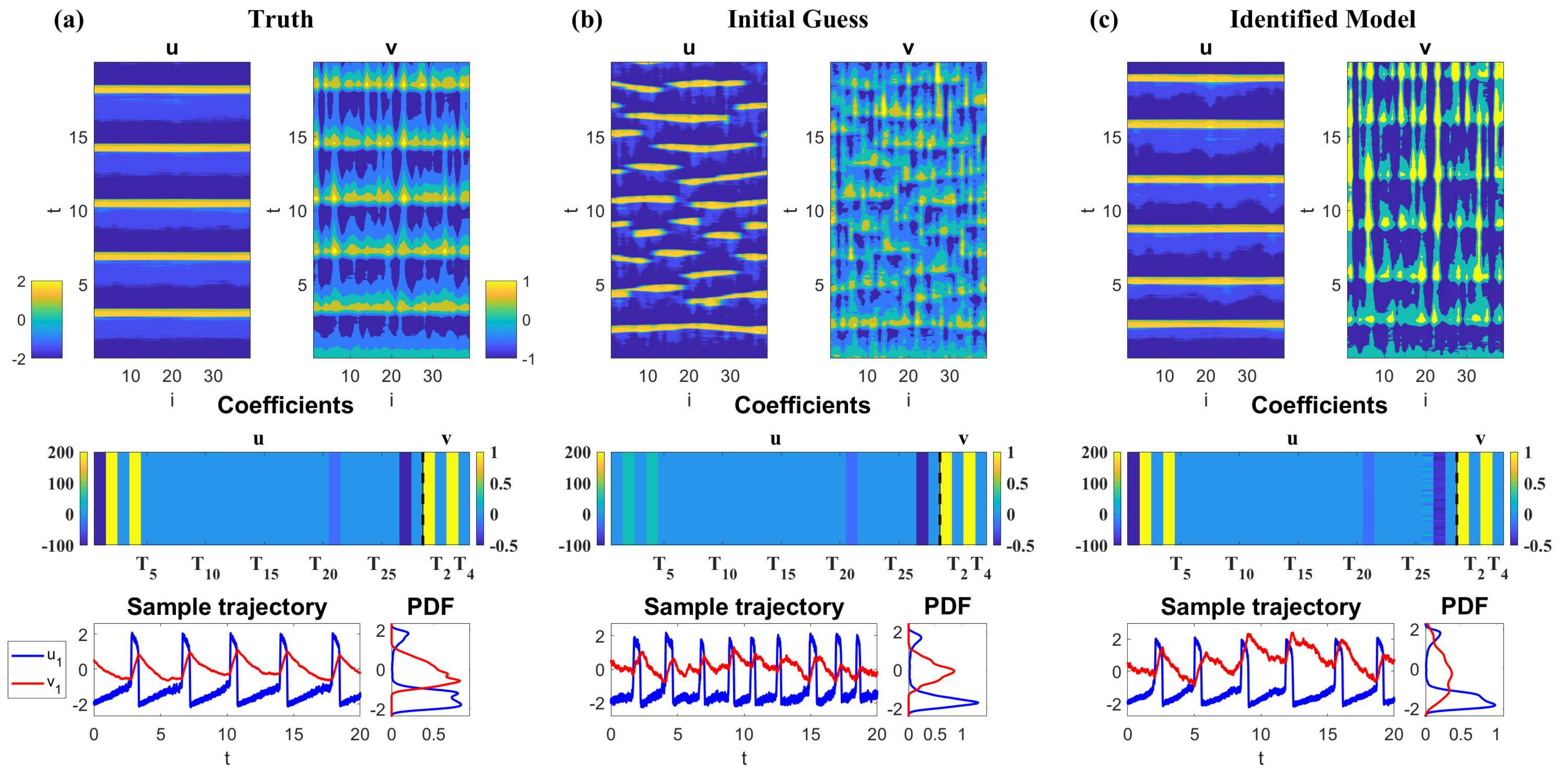}
\caption{Identification of the FHN model \eqref{FHN} for the first experiment. Panel (a): the true model structure. Panel (b): the initial guess. Panel (c): the identified model. In each panel, the first row shows the spatiotemporal patterns of the $u_i$ and the $v_i$. The second row shows the model structure and the associated parameters, where the $28$ columns of each $u_i$ correspond to the $28$ candidate functions in \eqref{FHN_Candidate_u} ordered in the same way and similar for the $4$ columns of each $v_i$. The third row shows the time series of $u_1$ and $v_1$.}
\label{FHN_comparison}
\end{figure}

\begin{figure}[!htb]
\centering
\includegraphics[width=1\textwidth]{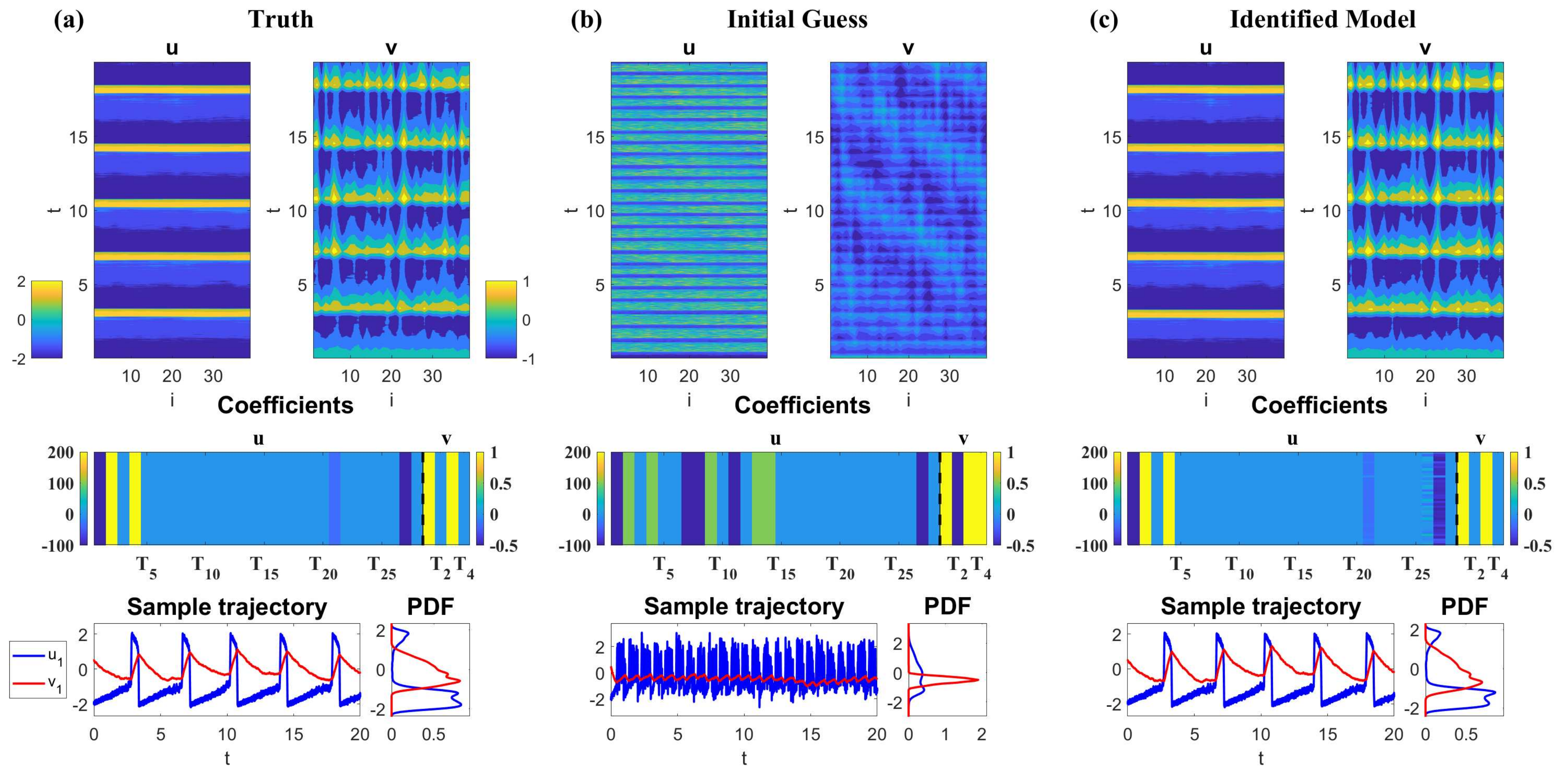}
\caption{Similar as Figure \ref{FHN_advection} but for the second experiment.}
\label{FHN_advection}
\end{figure}

The results of the two experiments are shown in Figure \ref{FHN_comparison} and Figure \ref{FHN_advection}, respectively, where Panel (c) in each figure displays the identified model.
In both experiments, the algorithm converges after $10$ iterations. The second row depicts the structures and the associated parameter values of the perfect model, the initial guess, and the identified model, where the $28$ columns of each $u_i$ correspond to the $28$ candidate functions in \eqref{FHN_Candidate_u} ordered in the same way and similar for the $4$ columns of each $v_i$.

From  Figure \ref{FHN_comparison}, it is seen that although the noise coefficient $\delta_2 = 0.4$ in the unobserved process is fixed and is chosen to be different from the truth $\delta_2 = 0.1$ that leads to a distinct spatiotemporal pattern of the $v_i$, the spatiotemporal pattern of the $u_i$ generated from the identified model is almost the same as the truth. Specifically, the identified model recovers the strong coherent spatiotemporal structure, which is completely missed in the initial guess of the model. The parameters in the identified model also have very similar values to the truth. In particular, the large error in the initial guess of the diffusion coefficient $d_u$, corresponding to coefficients  \#1, \#2, and \#4 in the figure (see also \eqref{FHN_Candidate_u}), is almost completely eliminated in the identified model. This is the key mechanism that leads to a strong coherent structure. Note that, as the $v_i$s are not observed, they are treated as the stochastic parameterizations in the identified model. Therefore, the $v_i$s in the identified model are not necessarily the same as the truth, but their statistical feedback to each $u_i$ is the crucial component that leads to the correct pattern of the latter, which the identified model is captured.

Figure \ref{FHN_advection} illustrates similar results. Despite the large difference in the model structure and the model parameters in the initial guess, the identified model recovers the truth accurately. Note that, as the true noise coefficient $\delta_2$ is utilized in the identified model, the pattern of $v_i$s is also clearly identified in addition to the recovery of the strong coherent structure of the $u_i$.

\section{Conclusions and Discussion}\label{Sec:Conclusions}
In this paper, a causality-based learning algorithm is developed that alternatives between model structure identification, conditional sampling of unobserved state variables, and parameter estimation to learn the underlying dynamics of the resolved scale variables as well as help develop suitable stochastic parameterization for the unresolved variables. Different from the constrained optimization with an L1 regularizer, the method developed here exploits the causation entropy to pre-determine the candidate functions that have potential contributions to the dynamics, which retains a quadratic optimization problem for parameter estimation via maximum likelihood estimates. The closed analytic formula of the conditional sampling allows an efficient recovery of the unobserved trajectories that facilitates the calculation of the causation entropy for the time evolutions of both the observed and the unobserved state variables. Physics constraints and localization techniques are further incorporated into the learning algorithm to include the basic physical properties in the data-driven models and reduce the computational cost, respectively. 

A hierarchy of chaotic and turbulent systems is adopted as test models. It is shown that the new learning algorithm effectively reproduces the dynamical and statistical features of the observed variables and provides suitable stochastic parameterizations with parsimonious structures. Many related topics are studied when implementing the numerical tests, including showing the necessity of the stochastic parameterization, understanding the effect of choosing different thresholds of the causation matrix selection, and detecting the importance of various terms in the original systems.

In addition to the library of the candidate functions, the work here does not assume any known model structure based on the prior knowledge. One natural extension of the current work is to learn the statistical closure of a turbulent system, where part of the model information is given. In addition, the iterative algorithm designed here only aims at finding a local optimum. This is sufficient for many applications. Yet, if the initial guess is very far from the truth or suitable local optimums, then the identified model may not fully capture the dynamical and statistical features of nature. Therefore, additional criteria can be incorporated into the learning procedure to facilitate the convergence of the algorithm at least towards a suitable local optimum.

\section*{Acknowledgements}
The research of N.C. is partially funded by ONR N00014-21-1-2904 and NSF 2118399. Y.Z. is supported as a PhD research assistant under these grants.



\bibliography{references}

\end{document}